\documentclass[11pt,oneside,reqno]{amsart}

\usepackage[margin=1.15in]{geometry}
\usepackage{amsmath,amssymb,amsthm,mathtools,xurl}
\usepackage{microtype}
\usepackage{enumitem}
\usepackage[colorlinks=true,linkcolor=blue,citecolor=blue,urlcolor=blue]{hyperref}

\newtheorem{theorem}{Theorem}[section]
\newtheorem{proposition}[theorem]{Proposition}
\newtheorem{lemma}[theorem]{Lemma}
\newtheorem{corollary}[theorem]{Corollary}
\theoremstyle{remark}
\newtheorem{remark}[theorem]{Remark}

\numberwithin{equation}{section}

\newcommand{\R}{\mathbb R}

\newcommand{\supp}{\operatorname{supp}}
\newcommand{\dist}{\operatorname{dist}}
\newcommand{\1}{\mathbf 1}

\title[Finite-Time Blow-up for Mass-Critical Fractional NLS]{Finite-Time Blow-up for the Radial Mass-Critical Fractional Nonlinear Schr\"odinger Equation}

\author[N.~U.~C.~Nguyen]{Ngoc Uyen Cong Nguyen}
\address[N.~U.~C. Nguyen]{Claire and Leo LLC, Houston, Texas, USA}
\email{contact@uyencong.com}

\author[V.~D.~Dinh]{Van Duong Dinh}
\address[V.~D.~Dinh]{Claire and Leo LLC, Houston, Texas, USA}
\email{contact@duongdinh.com}

\date{November 16, 2026}


\begin{document}

\begin{abstract}
We prove that radial negative-energy strong solutions to the focusing mass-critical fractional nonlinear Schr\"odinger equation
\[
i\partial_tu=(-\Delta)^su-|u|^{4s/N}u, \quad N\ge2,\quad \frac12<s<1,
\]
blow up in finite time. This resolves the finite-time blow-up versus infinite-time grow-up alternative left open by Boulenger, Himmelsbach and Lenzmann in the mass-critical regime. The argument uses an exponentially localized virial field together with a dyadic resolvent estimate that transfers the resulting weighted Balakrishnan energy to localized fractional Sobolev norms. When $s\ge N/(2N-1)$, we extend the conclusion to radial $H^s$ solutions by establishing the integrated virial identity at energy regularity. As a consequence, mass-critical ground-state standing waves are strongly unstable by finite-time blow-up.
\end{abstract}

\maketitle

\section{Introduction}
\label{sec:intro}

We consider the focusing mass-critical fractional nonlinear Schr\"odinger
equation
\begin{equation}
\left\{
\begin{aligned}
i\partial_tu&=(-\Delta)^su-|u|^{\frac{4s}{N}}u,\\
u(0)&=u_0\in H^s(\R^N),
\end{aligned}
\right.
\label{eq:fnls}
\end{equation}
where
\begin{equation}
N\ge2, \quad \frac12<s<1.
\label{eq:range}
\end{equation}

Throughout the paper, we write
\begin{equation} \label{eq:p-q}
p:=\frac{4s}{N},
\quad
q:=p+2=2+\frac{4s}{N}.
\end{equation}
Equation \eqref{eq:fnls} is mass-critical since it is invariant under the mass-preserving scaling
\[
u_\lambda(t,x)=\lambda^{N/2}u(\lambda^{2s}t,\lambda x).
\]
Along the flow, the mass
\[
M(u):=\|u\|_{L^2}^2
\]
and the energy
\[
E(u):=\frac12\|(-\Delta)^{s/2}u\|_{L^2}^2-\frac1q\|u\|_{L^q}^q
\]
are conserved.

The fractional Schr\"odinger equation was introduced by Laskin \cite{Laskin} in the framework of fractional quantum mechanics. Its
nonlinear dispersive theory has since developed in several directions. Local well-posedness in Sobolev spaces for power-type nonlinearities was studied by Hong and Sire \cite{HongSire} and the second author \cite{DinhWellposed}, among others. In the radial setting, the improved Strichartz estimates of Guo and Wang \cite{GuoWang} yield lossless estimates in a substantial range of fractional orders. Global well-posedness and scattering have also been investigated in critical and intercritical regimes; see, for example, Guo--Sire--Wang--Zhao \cite{GuoSireWangZhao} for the radial energy-critical equation and Sun--Wang--Yao--Zheng \cite{SunWangYaoZheng} for scattering below the ground-state threshold in the focusing case. On the elliptic side, the uniqueness and nondegeneracy of radial ground states for fractional equations were established by Frank, Lenzmann and Silvestre \cite{FrankLenzmannSilvestre} and related works.

Our concern here is finite-time blow-up. For the classical nonlinear Schr\"odinger equation $s=1$, virial and localized virial identities are among the fundamental tools for proving blow-up; see, for example, Glassey \cite{Glassey} and Ogawa--Tsutsumi
\cite{OgawaTsutsumi}. The fractional problem is substantially more delicate because $(-\Delta)^s$ is nonlocal. In particular, the usual
local virial calculus does not apply directly, and spatial localization creates interactions between different regions of space.

A fundamental advance was made by Boulenger, Himmelsbach and Lenzmann \cite{BHL}, who developed a localized virial method for fractional NLS based on the Balakrishnan representation of the fractional Laplacian. In the radial mass-critical case, they showed that negative energy implies the alternative
\[
T_+<\infty,
\]
or
\[
T_+=\infty \quad\text{and}\quad \|(-\Delta)^{s/2}u(t)\|_{L^2}\gtrsim t^s
\]
for all sufficiently large $t$. Thus their argument does not exclude a global solution whose fractional kinetic norm grows without bound.

The blow-up theory for fractional Schr\"odinger equations has since attracted a considerable interest in the mathematical community. For power-type nonlinearities, Feng \cite{FengCombined} studied equations with combined powers, obtaining blow-up criteria, sharp global-existence/blow-up thresholds, and concentration and profile properties of blow-up solutions. In \cite{DinhSupercritical,DinhMassCritical} the second author investigated blow-up dynamics in the mass-supercritical and mass-critical regimes. He also established localized-virial blow-up criteria beyond radial data in \cite{DinhCriteria}. The instability of standing waves was studied by Peng and Shi \cite{PengShi} in the mass-critical setting and by the second author \cite{DinhInstability} in the mass-supercritical case.

Related blow-up phenomena have also been studied for fractional Schr\"odinger equations with nonlocal nonlinearities. Guo and Zhu
\cite{GuoZhuHartree} established a sharp radial threshold separating scattering and finite-time blow-up for the $L^2$-supercritical
fractional Hartree equation. Feng, Chen and Liu \cite{FengChenLiuChoquard} obtained nonradial blow-up criteria and strong instability results for the fractional Schr\"odinger--Choquard equation. Earlier concentration and blow-up-profile results for mass-critical fractional Hartree equations were obtained, among others, by Cho--Hwang--Kwon--Lee \cite{ChoHwangKwonLeeProfile}.

Of particular relevance to the mass-critical problem, Lan \cite{LanBlowupDynamics} proved finite-time blow-up and obtained a
detailed description of the blow-up dynamics for the one-dimensional $L^2$-critical fractional Schr\"odinger equation when the fractional dispersion exponent is sufficiently close to the classical Schr\"odinger case. More recently, Gou--R\u{a}dulescu--Zhang
\cite{GouRadulescuZhang} extended localized-virial blow-up criteria for fractional NLS from radial to cylindrically symmetric data. Their mass-critical result, however, still retains the finite-time blow-up/infinite-time grow-up alternative. Thus, for the pure-power problem in dimensions \(N\ge2\), the global grow-up branch in the general radial negative-energy mass-critical theory remains unresolved.

The purpose of this paper is to rule out this grow-up scenario for radial negative-energy solutions. Our main result is the following.

\begin{theorem}[Finite-time blow-up for strong radial solutions]
\label{thm:strong}
Assume \eqref{eq:range}. Let
\[
u\in C((-T_-,T_+);H^{2s}_{\rm rad}(\R^N))
\]
be a maximal strong solution of \eqref{eq:fnls}. If $E(u_0)<0$, then
\[
T_+<\infty, \quad T_-<\infty.
\]
\end{theorem}

The $H^{2s}$ regularity in Theorem~\ref{thm:strong} is used only to justify the localized virial identity at the strong-solution level.
All estimates obtained from the virial identity involve quantities at the $H^s$ level. By establishing the integrated virial identity for $H^s$ solution, we obtain the following energy-space formulation.

\begin{corollary}[Energy-space formulation]
\label{cor:energy-space}
Assume
\begin{equation} \label{eq:expo-energy}
N\ge2, \quad \frac{N}{2N-1}\le s<1.
\end{equation}
Let $u_0\in H^s_{\rm rad}(\R^N)$, and let $u \in C((-T_-,T_+);H^{s}_{\rm rad}(\R^N))$ be the maximal solution furnished by the corresponding $H^s$ local well-posedness theory. If $E(u_0)<0$, then
\[
T_+<\infty, \quad T_-<\infty.
\]
In particular, the $H^s$ norm diverges at each finite endpoint of the maximal lifespan.
\end{corollary}

The restriction $N\ge2$ in Theorem~\ref{thm:strong} enters through the radial Sobolev estimate used to control the nonlinear localization error, which requires
\[
\frac12<s<\frac N2.
\]
The additional lower bound
\[
s\ge\frac{N}{2N-1}
\]
in Corollary~\ref{cor:energy-space} comes from the currently available lossless radial local theory rather than from the virial argument
itself; see \cite{GuoWang,DinhMassCritical}.

As a consequence, the known strong instability of the mass-critical ground-state standing waves can be strengthened to strong instability by genuine finite-time blow-up.

\begin{corollary}[Strong instability by finite-time blow-up]
\label{cor:instability}
Assume \eqref{eq:expo-energy}. Let $Q_\omega\in H^s_{\rm rad}(\R^N)$ be a positive ground state satisfying
\begin{equation}
(-\Delta)^sQ_\omega +\omega Q_\omega -Q_\omega^{1+\frac{4s}{N}} =0,
\quad
\omega>0.
\label{eq:ground-state}
\end{equation}
Then, for every $\varepsilon>0$, there exists $u_0\in H^s_{\rm rad}(\R^N)$ such that 
\[
\|u_0-Q_\omega\|_{H^s}<\varepsilon
\]
and the corresponding solution of \eqref{eq:fnls} blows up in finite time.
\end{corollary}

Let us describe the main ideas behind Theorem~\ref{thm:strong}. 

A finite-time mechanism that has recently emerged in several mass-critical dispersive problems is the use of an exponentially localized virial field. In \cite{NguyenDinhQuadratic}, the authors introduced such a multiplier for the four-dimensional mass-critical quadratic Schr\"odinger system without mass resonance. The resulting weighted kinetic defect, combined with a weighted interpolation estimate and a superlinear differential inequality, allows one to exclude the global grow-up scenario. Leonard Dinh \cite{LeonardDinhBiharmonic} subsequently employed a related exponential localization for the homogeneous mass-critical biharmonic NLS and proved finite-time blow-up for radial negative-energy solutions, thereby eliminating the mass-critical grow-up alternative left open in the earlier biharmonic theory of Boulenger and Lenzmann \cite{BL}.  More recently, in \cite{NguyenDinhInhomBiharmonic} the authors developed an analogous mechanism for the mass-critical inhomogeneous biharmonic NLS and obtained nonradial finite-time blow-up.

The essential new difficulty in the present problem is the nonlocality of the fractional dispersive operator. For second- and fourth-order local operators, the defect generated by a spatial multiplier can be expressed directly through local derivatives of the solution. For $(-\Delta)^s$, this is no longer possible. The Balakrishnan representation rewrites the fractional kinetic energy in terms of the resolvent family
\[
u_\mu=c_s(-\Delta+\mu)^{-1}u,
\]
and the localized fractional virial identity naturally produces a weighted resolvent quantity rather than a directly localized
$\dot H^s$ norm. Since the resolvent is itself nonlocal, multiplication by a spatial cutoff does not commute with it, and different dyadic regions continue to interact.

The mass-critical nonlinearity creates an additional difficulty. In the localized virial argument of Boulenger--Himmelsbach--Lenzmann \cite{BHL}, the localization function $\phi_R$ agrees with the quadratic function on a large ball. Consequently,
\[
I-D^2\phi_R=0
\]
throughout this region, so the localization does not provide any control of the kinetic energy there. This is not a serious obstacle in the mass-supercritical case, where the localized virial identity contains a favorable coercive term even without such control. At the mass-critical exponent, however, the leading kinetic and nonlinear contributions cancel exactly. As a result, the localized virial estimate still shows that the virial observable decreases, but it does not provide sufficient control of the observable in terms of the remaining kinetic contribution. This is why the usual argument yields a blow-up-or-grow-up alternative rather than ruling out global solutions altogether.

Motivated by the exponential-localization mechanism above, we introduce
\[
\nabla\phi_R(x)=x\exp\left[-\left(\frac{|x|}{R}\right)^{2s}\right]
\]
and
\[
\omega_R(x) = 1-\exp\left[-\left(\frac{|x|}{R}\right)^{2s}\right].
\]
The fractional virial identity naturally controls the Balakrishnan defect
\[
\mathcal Q_R(u):= \int_0^\infty \mu^s \int_{\R^N} \omega_R(x)|\nabla u_\mu(x)|^2 dx d\mu.
\]

To describe its localized Sobolev counterpart, let $\rho_j=2^jR$, choose radial annular cutoffs $\zeta_j,\eta_j$ such that
\[
\sum_{j\in\mathbb Z}\zeta_j^2=1, \quad \eta_j\equiv1\quad\text{on }\supp\zeta_j,
\]
and set
\[
d_j:=\min\left\{\left(\frac{\rho_j}{R}\right)^{2s},1\right\}.
\]
Formally, the dyadic fractional defect is
\[
\mathcal D_R(u):= \sum_{j\in\mathbb Z}d_j \left(\|(-\Delta)^{s/2}(\zeta_j u)\|_{L^2}^2+\|(-\Delta)^{s/2}(\eta_j u)\|_{L^2}^2\right).
\]
The precise construction of the cutoff families is given in Section~\ref{sec:defect}.

The first main step is the Balakrishnan-to-dyadic transfer estimate
\begin{equation}\label{eq:intro-transfer}
\mathcal D_R(u) \lesssim \mathcal Q_R(u)+R^{-2s}M(u).
\end{equation}
It is proved by combining a dyadic decomposition with off-diagonal resolvent decay between separated annuli. This transfer is the principal new ingredient of the present argument: it converts the nonlocal coercivity furnished by the
fractional virial identity into localized fractional Sobolev control.

The second step is to control the nonlinear localization error. Radial Sobolev estimates on the dyadic annuli yield, for every
$\varepsilon>0$,
\begin{equation}
\int_{\R^N} \omega_R(x)|u(x)|^q dx \le \varepsilon \mathcal D_R(u) + C_\varepsilon R^{-2s} M(u)^{\frac{N}{N-2s}}.
\label{eq:intro-nonlinear}
\end{equation}
Combining \eqref{eq:intro-transfer} and \eqref{eq:intro-nonlinear} with the fractional virial identity gives the coercive localized virial estimate
\[
\mathcal M_R'(t) \le 8sE(u_0) -c\mathcal D_R(u(t)) + CR^{-2s} \left[M(u_0)+M(u_0)^{\frac{N}{N-2s}}\right],
\]
where
\[
\mathcal M_R(u):=2 \operatorname{Im} \int_{\R^N} \overline u \nabla\phi_R\cdot\nabla u dx.
\]
If $E(u_0)<0$, we may fix $R$ sufficiently large so that
\[
\mathcal M_R'(t) \le -c-c\mathcal D_R(u(t)).
\]

The final step couples the same defect back to the virial observable. We prove
\begin{equation}
|\mathcal M_R(u)| \lesssim M(u) + R M(u)^{1-\frac1{2s}} \mathcal D_R(u)^{\frac1{2s}}.
\label{eq:intro-observable}
\end{equation}
Once $\mathcal M_R(t)$ becomes sufficiently negative, \eqref{eq:intro-observable} yields a lower bound for $\mathcal D_R(u(t))$ in terms of $-\mathcal M_R(t)$. Setting
\[
y(t):=-\mathcal M_R(u(t)),
\]
we obtain
\[
y'(t)\gtrsim y(t)^{2s}.
\]
Since $2s>1$, this superlinear Osgood inequality cannot persist for all positive times. Hence the global grow-up branch is impossible and the solution must blow up in finite time.

\begin{remark}
\label{rem:supercritical}
The exponential-defect mechanism is not specific to the mass-critical nonlinearity. In the mass-supercritical and energy-critical regimes, the localized virial identity already retains a favorable coercive contribution after combining the leading kinetic and nonlinear terms. The exponential multiplier above is compatible with this coercivity, and the dyadic Balakrishnan localization controls the additional nonlocal localization errors. Thus the same strategy is compatible with these regimes and appears to provide an alternative route to the corresponding finite-time blow-up results of Boulenger--Himmelsbach--Lenzmann \cite{BHL}. We do not pursue these cases here, since the distinctive difficulty of the present problem occurs at the mass-critical exponent, where the global coercive contribution degenerates and the defect must additionally control the virial observable itself.
\end{remark}

The paper is organized as follows. Section~\ref{sec:prelim} collects the local theory, the Balakrishnan representation, the fractional virial identity, and the analytic inequalities used throughout the paper. In Section~\ref{sec:weight}, we introduce the exponential multiplier, establish its geometric properties, and justify the regularization needed in the fractional virial identity. Section~\ref{sec:defect} develops the dyadic fractional defect, proves the off-diagonal resolvent estimate and the Balakrishnan-to-dyadic defect transfer, and establishes the radial weighted interpolation estimate. Section~\ref{sec:virial} derives the coercive localized fractional virial inequality. Section~\ref{sec:virial-defect} establishes the virial--defect estimate.
The finite-time blow-up theorem is proved in Section~\ref{sec:blowup}, and Section~\ref{sec:instability} is devoted to strong instability by finite-time blow-up of the ground-state standing waves.

\section{Preliminaries}
\label{sec:prelim}

Throughout the paper,
\[
H^s_{\mathrm{rad}}(\mathbb R^N) :=\left\{f\in H^s(\mathbb R^N): f(x)=F(|x|) \text{ for some function }F\right\}
\]
denotes the closed subspace of radial functions in \(H^s(\mathbb R^N)\).

\subsection{Local theory}

The local well-posedness theory for fractional nonlinear Schr\"odinger equations in Sobolev spaces was developed, among others, by Hong and Sire \cite{HongSire} and the second author \cite{DinhWellposed}. For radial data, the improved Strichartz estimates of Guo and Wang \cite{GuoWang} yield a lossless $H^s$ theory in the range $s\ge N/(2N-1)$; see \cite[Proposition 2.5]{DinhMassCritical}.

\begin{proposition}[Local theory in the energy space \cite{DinhMassCritical}]
\label{prop:local-theory}
Assume \eqref{eq:expo-energy}. For every radial $u_0\in H^s(\R^N)$ there exist maximal times $T_\pm\in(0,\infty]$ and a unique solution 
\[
u\in C((-T_-,T_+);H^s_{\rm rad}(\R^N)) \cap L^\gamma_{\rm loc}((-T_-,T_+);W^{s,\rho}(\R^N)),
\]
to \eqref{eq:fnls}, where
\begin{equation} \label{eq:gamma-rho}
\gamma=\frac{2(N+2s)}{N-2s}, \quad \rho=\frac{2N(N+2s)}{N^2+4s^2}.
\end{equation}
Moreover, mass and energy are conserved, and the usual blow-up alternative holds: if $T_+<\infty$, then
\[
\lim_{t\uparrow T_+}\|u(t)\|_{H^s}=\infty,
\]
with the analogous statement backward in time.
\end{proposition}

\subsection{Balakrishnan representation and the fractional virial identity}

Let
\begin{equation}
u_\mu:=c_s(-\Delta+\mu)^{-1}u,
\quad
c_s:=\left(\frac{\sin(\pi s)}{\pi}\right)^{1/2}.
\label{eq:umu}
\end{equation}
Balakrishnan's representation \cite{Balakrishnan} and Plancherel's theorem imply (see \cite[Equation (2.12)]{BHL})
\begin{equation} \label{eq:balakrishnan-energy}
\int_0^\infty \mu^s\|\nabla u_\mu\|_{L^2}^2 d\mu=s\|(-\Delta)^{s/2}u\|_{L^2}^2.
\end{equation}
We shall repeatedly use the elementary resolvent bounds
\begin{equation}
\|(-\Delta+\mu)^{-1}\|_{L^2\to L^2}\le \mu^{-1},
\quad
\|\nabla(-\Delta+\mu)^{-1}\|_{L^2\to L^2}\lesssim \mu^{-1/2}.
\label{eq:resolvent-bounds}
\end{equation}

For a real-valued weight $\phi$, define the localized virial action
\[
\mathcal M_\phi(u):=2\operatorname{Im}\int_{\R^N}\overline{u} \nabla\phi\cdot\nabla u dx.
\]
For a real-valued function $\phi$ with $\nabla \phi \in W^{3,\infty}(\mathbb R^N)$, Boulenger--Himmelsbach--Lenzmann \cite[Lemma 2.1]{BHL} established the following localized virial identity for $i\partial_tu=(-\Delta)^su-|u|^p u$ under \eqref{eq:range}: 
\begin{align}
\frac{d}{dt} \mathcal M_\phi(u(t))
={}&4\int_0^\infty\mu^s\int_{\R^N}
\partial_k\overline{u_\mu} \partial_{kl}\phi \partial_lu_\mu dx d\mu \notag\\
&-\int_0^\infty\mu^s\int_{\R^N}(\Delta^2\phi)|u_\mu|^2 dx d\mu
-\frac{2p}{p+2}\int_{\R^N}(\Delta\phi)|u|^{p+2} dx.
\label{eq:BHL-virial}
\end{align}
Repeated indices in \eqref{eq:BHL-virial} are summed.

The identity above is first justified for $H^{2s}$-valued strong solutions, as in \cite{BHL}. A further justification is needed for our exponential localization function $\phi_R$, since its fourth-order derivatives are singular at the origin and therefore do not satisfy the regularity assumptions of the identity as stated above. To overcome this difficulty, we first apply \eqref{eq:BHL-virial} to the smooth approximation $\phi_{R,\delta}$ introduced in Section~\ref{subsec:regularization}, and then pass to the limit as $\delta\downarrow0$. This procedure allows us to derive the localized virial estimates for strong solutions without applying the original identity directly to the singular multiplier. We subsequently extend the argument to radial solutions in the energy
space $H^s(\R^N)$ by establishing the corresponding integrated virial identity in Appendix~\ref{app:energy-approx}.

\subsection{Radial Sobolev and Hardy inequalities}

We use the Cho--Ozawa radial Sobolev inequality \cite{ChoOzawa}: if
$f\in\dot H^s_{\mathrm{rad}}(\R^N)$ and
\[
\frac12<s<\frac N2,
\]
then
\begin{equation}
\sup_{x\ne0}|x|^{N/2-s}|f(x)|
\le C_{N,s}\|(-\Delta)^{s/2}f\|_{L^2}.
\label{eq:radial-sobolev}
\end{equation}
Under \eqref{eq:range}, the condition $s<N/2$ is automatic.

We also use the fractional Hardy inequality \cite{FrankSeiringer}: for
$0<\alpha<N/2$,
\begin{equation}
\||x|^{-\alpha}f\|_{L^2}
\lesssim \|(-\Delta)^{\alpha/2}f\|_{L^2}.
\label{eq:hardy}
\end{equation}

\section{The exponential localization and its geometry}
\label{sec:weight}

For $R>0$, let $r=|x|$, $z=r/R$, and set
\begin{equation}
h(z):=e^{-z^{2s}}, \quad \phi_R'(r):=rh(r/R).
\label{eq:weight-def}
\end{equation}
Thus
\[
\nabla\phi_R(x)=xh(|x|/R), \quad \omega_R(x):=1-h(|x|/R).
\]
The elementary inequality $1-e^{-y}\simeq\min\{y,1\}$ for $y\ge0$ gives
\begin{equation}
\omega_R(x)\simeq
\min\left\{\left(\frac{|x|}{R}\right)^{2s},1\right\}.
\label{eq:defect-size}
\end{equation}

\subsection{Geometry of the multiplier}

\begin{lemma}[Exponential defect geometry]
\label{lem:geometry}
Let $N\geq 1$ and $0<s<1$. There exists $C=C(N,s)$ such that for every $R>0$,
\begin{align}
I-D^2\phi_R&\ge \omega_R I,
\label{eq:hessian-defect}\\
0\le N-\Delta\phi_R&\le C\omega_R,
\label{eq:lap-defect}\\
|\nabla\phi_R(x)| &\le CR \left(\omega_R(x)\right)^{\frac{1}{2s}},
\label{eq:gradient-defect}\\
|\Delta^2\phi_R(x)|&\le CR^{-2s}|x|^{-2(1-s)},
\quad x\ne0.
\label{eq:bilap-weight}
\end{align}
Moreover,
\[
\|\nabla\phi_R\|_{L^\infty}\lesssim R,
\quad
\|D^2\phi_R\|_{L^\infty}\lesssim1.
\]
\end{lemma}

\begin{proof}
For a radial function, the tangential and radial eigenvalues of the Hessian are $\phi_R'(r)/r$ and $\phi_R''(r)$. From \eqref{eq:weight-def},
\[
\frac{\phi_R'(r)}r=h(z),
\quad
\phi_R''(r)=h(z)(1-2sz^{2s}).
\]
Hence the tangential eigenvalue of $I-D^2\phi_R$ is $1-h=\omega_R$, whereas the radial eigenvalue is
\[
1-h(1-2sz^{2s})=\omega_R+2sz^{2s}h\ge\omega_R.
\]
This proves \eqref{eq:hessian-defect}. Also
\[
\Delta\phi_R= \phi''_R + (N-1)\frac{\phi'_R}{r}=h(z)(N-2sz^{2s}),
\]
so
\[
N-\Delta\phi_R=N\omega_R+2sz^{2s}e^{-z^{2s}}.
\]
Since $ye^{-y}\lesssim1-e^{-y}$ for $y\ge0$, \eqref{eq:lap-defect} follows.

If $r\le R$, then \eqref{eq:defect-size} gives $\omega_R\simeq(r/R)^{2s}$, hence
\[
|\nabla\phi_R| \omega_R^{-1/(2s)} \simeq R \frac{|x|}{R} e^{-\left(\frac{|x|}{R}\right)^{2s}} \left(\frac{|x|}{R}\right)^{-1} \lesssim R.
\]
If $r\ge R$, then $\omega_R\simeq1$ and $re^{-(r/R)^{2s}}\lesssim R$. This proves \eqref{eq:gradient-defect}.

To prove \eqref{eq:bilap-weight}, write $\phi_R(x)=R^2\Phi(x/R)$ with $\Phi$ is radial satisfying
\[
\Phi'(r)=rH(r), \quad H(r)=e^{-r^{2s}}.
\]
For a radial function with $\Phi'(r)=rH(r)$, direct differentiation gives
\[
\Delta^2\Phi(r)=rH'''(r)+(2N+1)H''(r) +\frac{N^2-1}{r}H'(r).
\]
Substituting $H(r)=e^{-r^{2s}}$ yields
\begin{align*}
\Delta^2\Phi(r)
=e^{-r^{2s}}r^{2s-2}\bigl[-2s(N+2s)(N+2s-2) +(2s)^2(2N+6s-2)r^{2s}
-(2s)^3r^{4s}\bigr].
\end{align*}
The polynomial factor is bounded after multiplication by the exponential, and hence
\[
|\Delta^2\Phi(r)|\le C_{N,s}r^{2s-2},\quad r>0.
\]
By scaling,
\[
|\Delta^2\phi_R(x)| =R^{-2}|\Delta^2\Phi(|x|/R)| \le C_{N,s}R^{-2s}|x|^{2s-2}.
\]
This proves \eqref{eq:bilap-weight}. The remaining bounds follow directly from \eqref{eq:weight-def}.
\end{proof}

As $\frac12<s<1$, the fourth-order derivatives of $\phi_R$ are singular at the origin. The estimate \eqref{eq:bilap-weight} puts this singularity in the range of the fractional Hardy inequality. We therefore regularize the multiplier before applying the virial identity.

\subsection{Regularization of the multiplier}
\label{subsec:regularization}

For $0<\delta<R$, define
\[
g_{R,\delta}(r):=R^{-2s}\big((r^2+\delta^2)^s-\delta^{2s}\big),
\quad
h_{R,\delta}(r):=e^{-g_{R,\delta}(r)},
\]
and let
\[
\phi_{R,\delta}'(r):=r h_{R,\delta}(r).
\]
Then $\phi_{R,\delta}$ is smooth and radial, and
\[
\nabla\phi_{R,\delta}\longrightarrow\nabla\phi_R,
\quad
D^2\phi_{R,\delta}\longrightarrow D^2\phi_R
\]
locally uniformly away from the origin and pointwise everywhere. Moreover, with constants independent of $\delta$,
\[
\|\nabla\phi_{R,\delta}\|_{L^\infty}\lesssim R,
\quad
\|D^2\phi_{R,\delta}\|_{L^\infty}\lesssim1,
\]
and
\begin{equation}
|\Delta^2\phi_{R,\delta}(x)|
\lesssim R^{-2s}(|x|+\delta)^{-2(1-s)}.
\label{eq:regularized-bounds2}
\end{equation}
To justify \eqref{eq:regularized-bounds2}, set
\[
\varepsilon:=\frac{\delta}{R}, \quad y:=\frac{x}{R}.
\]
Then, up to an irrelevant additive constant,
\[
\phi_{R,\delta}(x)=R^2\Phi_\varepsilon(y),
\]
where $\Phi_\varepsilon$ is radial and
\[
\Phi_\varepsilon'(r)= r\exp \left(-\bigl[(r^2+\varepsilon^2)^s-\varepsilon^{2s}\bigr]\right).
\]
Hence
\[
\Delta_x^2\phi_{R,\delta}(x)=R^{-2}\Delta_y^2\Phi_\varepsilon(y).
\]
Writing
\[
G_\varepsilon(r):=(r^2+\varepsilon^2)^s-\varepsilon^{2s},
\]
direct differentiation gives, for $k=1,2,3$,
\[
|G_\varepsilon^{(k)}(r)| \lesssim (r+\varepsilon)^{2s-k}.
\]
In the radial identity
\[
\Delta^2\Phi_\varepsilon
=rH_\varepsilon'''+(2N+1)H_\varepsilon''
+\frac{N^2-1}{r}H_\varepsilon',
\quad H_\varepsilon=e^{-G_\varepsilon},
\]
the quotient in the last term is controlled separately. Indeed, for $r>0$,
\[
\frac{|G_\varepsilon'(r)|}{r}
=2s(r^2+\varepsilon^2)^{s-1}
\lesssim(r+\varepsilon)^{2s-2}.
\]
On $0<r\le1$, these bounds, together with
\[
H_\varepsilon'=-G_\varepsilon'e^{-G_\varepsilon},\quad
H_\varepsilon''=\big((G_\varepsilon')^2-G_\varepsilon''\big)e^{-G_\varepsilon}
\]
(with the analogous formula for $H_\varepsilon'''$) therefore give
\[
|\Delta^2\Phi_\varepsilon(r)|\lesssim(r+\varepsilon)^{2s-2}.
\]
For $r\ge1$, the lower bound $G_\varepsilon(r)\ge r^{2s}-1$ and the uniform polynomial bounds on its first three derivatives show that $e^{-G_\varepsilon(r)}$ dominates the derivative factors. Combining the two regions yields
\[
|\Delta^2\Phi_\varepsilon(y)|\lesssim(|y|+\varepsilon)^{2s-2}
\quad(0<\varepsilon\le1).
\]
Therefore
\[
|\Delta^2\phi_{R,\delta}(x)| \lesssim R^{-2} \left(\frac{|x|+\delta}{R}\right)^{2s-2} =R^{-2s}(|x|+\delta)^{-2(1-s)},
\]
which proves \eqref{eq:regularized-bounds2}.

The tangential and radial eigenvalues of $D^2\phi_{R,\delta}$ are respectively
\[
h_{R,\delta}(r),
\quad
h_{R,\delta}(r)\big(1-rg_{R,\delta}'(r)\big),
\]
so that
\[
I-D^2\phi_{R,\delta}\ge (1-h_{R,\delta})I.
\]

We now justify the use of the exponential weight $\phi_R$ in the differential virial identity for strong solutions. Although
$\phi_R$ has bounded first- and second-order derivatives, its fourth-order derivatives are singular at the origin, so it does
not satisfy the regularity assumptions of \eqref{eq:BHL-virial}. We therefore first apply the identity to the smooth approximation $\phi_{R,\delta}$ introduced in Section~\ref{subsec:regularization}, and then remove the regularization.

Fix a compact interval $I$ contained in the lifespan of a strong solution $u\in C(I;H^{2s}(\R^N))$, and let
$t_1<t_2$ belong to $I$. For every $\delta>0$, the multiplier $\phi_{R,\delta}$ satisfies the hypotheses of
\eqref{eq:BHL-virial}. Integrating that identity in time, we obtain 
\[
\mathcal M_{\phi_{R,\delta}}(u(t_2))
-\mathcal M_{\phi_{R,\delta}}(u(t_1))
=
\int_{t_1}^{t_2}
\mathcal V_{R,\delta}(u(t)) dt,
\]
where $\mathcal V_{R,\delta}$ denotes the right-hand side of \eqref{eq:BHL-virial} with $\phi=\phi_{R,\delta}$. We shall
pass to the limit $\delta\downarrow0$ in this integrated identity before recovering the corresponding differential
identity.

We begin with the principal Balakrishnan term. The regularized Hessians converge pointwise to $D^2\phi_R$ and satisfy
\[
\sup_{\delta>0}
\|D^2\phi_{R,\delta}\|_{L^\infty}<\infty.
\]
Consequently, the absolute value of its integrand is bounded by $C\mu^s|\nabla u_\mu|^2$, independently of $\delta$.
The Balakrishnan identity \eqref{eq:balakrishnan-energy} gives 
\[
\int_0^\infty
\mu^s\|\nabla u_\mu(t)\|_{L^2}^2 d\mu
=
s\|(-\Delta)^{s/2}u(t)\|_{L^2}^2.
\]
Since $u\in C(I;H^{2s}(\R^N))$, the right-hand side is uniformly bounded on $I$. Thus the dominating function is integrable
in time, space, and the resolvent parameter $\mu$, and dominated convergence yields the desired limit of the
principal term. The nonlinear contribution is treated similarly, since $\Delta\phi_{R,\delta}$ converges pointwise
to $\Delta\phi_R$, remains uniformly bounded, and 
\[
\int_I\|u(t)\|_{L^q}^q dt<\infty
\]
by the Sobolev embedding $H^s(\R^N)\hookrightarrow L^q(\R^N)$.

The lower-order Balakrishnan term requires additional care because the fourth-order derivatives of the regularized
weights are not uniformly bounded as $\delta\downarrow0$. Instead, \eqref{eq:regularized-bounds2} provides the estimate
\[
|\Delta^2\phi_{R,\delta}(x)|
\lesssim R^{-2s}|x|^{-2(1-s)},
\quad x\neq0,
\]
with a constant independent of $\delta$. To justify dominated convergence, it is therefore sufficient to
establish the integrability of 
\[
\int_I\int_0^\infty\mu^s
\int_{\R^N}
|x|^{-2(1-s)}|u_\mu(t,x)|^2
 dx d\mu dt.
\]
Since $0<1-s<N/2$, the fractional Hardy inequality gives
\[
\int_{\R^N}
|x|^{-2(1-s)}|u_\mu|^2 dx
\lesssim
\|(-\Delta)^{(1-s)/2}u_\mu\|_{L^2}^2.
\]
Using the definition $u_\mu=c_s(-\Delta+\mu)^{-1}u$ and Plancherel's theorem, we obtain
\begin{equation} \label{eq:iden-L2}
\begin{aligned}
&\int_0^\infty
\mu^s\|(-\Delta)^{(1-s)/2}u_\mu\|_{L^2}^2 d\mu\\
&\quad=
c_s^2\int_{\R^N}
|\widehat u(\xi)|^2|\xi|^{2(1-s)}
\left(
\int_0^\infty
\frac{\mu^s}{(|\xi|^2+\mu)^2} d\mu
\right)d\xi\\
&\quad=
c_s^2
\left(
\int_0^\infty\frac{\tau^s}{(1+\tau)^2} d\tau
\right)\|u\|_{L^2}^2
=s\|u\|_{L^2}^2.
\end{aligned}
\end{equation}
Combining these estimates and integrating in time, we deduce that
\[
\int_I\int_0^\infty\mu^s
\int_{\R^N}
|x|^{-2(1-s)}|u_\mu(t,x)|^2
 dx d\mu dt
\lesssim
|I|\sup_{t\in I}\|u(t)\|_{L^2}^2<\infty.
\]
We have thus obtained an integrable majorant for the lower-order term that is independent of $\delta$. Since $\Delta^2\phi_{R,\delta}$ converges pointwise away from the origin to $\Delta^2\phi_R$, dominated convergence applies to this term as well.

It remains to justify convergence of the virial observables at the endpoints $t_1$ and $t_2$. The vector fields
$\nabla\phi_{R,\delta}$ converge pointwise to $\nabla\phi_R$ and are uniformly bounded in $W^{1,\infty}(\R^N)$ for fixed $R$. For Schwartz functions, convergence of the corresponding virial forms follows directly by dominated convergence. Lemma~\ref{lem:first-order-form} shows that these forms are uniformly bounded quadratic forms on $H^{1/2}(\R^N)$. We may therefore approximate
an arbitrary $H^s$-function by Schwartz functions and use $H^s(\R^N)\hookrightarrow H^{1/2}(\R^N)$ to conclude that
\[
\mathcal M_{\phi_{R,\delta}}(u(t_j))
\longrightarrow
\mathcal M_R(u(t_j)),
\quad j=1,2.
\]
Passing to the limit in the integrated identity, we obtain
\[
\mathcal M_R(u(t_2))-\mathcal M_R(u(t_1))
=
\int_{t_1}^{t_2}\mathcal V_R(u(t)) dt,
\]
where $\mathcal V_R$ denotes the right-hand side of \eqref{eq:BHL-virial} with the limiting weight $\phi_R$.

We finally show that this integrated identity implies the differential identity needed below. For this purpose, it
suffices to establish that $\mathcal V_R$ is continuous on $H^s(\R^N)$. The principal Balakrishnan contribution
defines a continuous quadratic form on $H^s(\R^N)$, because $D^2\phi_R$ is bounded and \eqref{eq:balakrishnan-energy} controls its associated bilinear form. Similarly, the Hardy--Balakrishnan estimate proved above shows that the lower-order contribution
defines a bounded quadratic form on $L^2(\R^N)$, whose continuity follows by polarization. Finally, the nonlinear functional
\[
v\longmapsto
\int_{\R^N}(\Delta\phi_R)|v|^q dx
\]
is continuous on $H^s(\R^N)$, since $\Delta\phi_R\in L^\infty(\R^N)$ and $H^s(\R^N)\hookrightarrow L^q(\R^N)$. It follows that $\mathcal V_R:H^s(\R^N)\to\R$ is continuous.

Since $u\in C(I;H^{2s}(\R^N))\subset C(I;H^s(\R^N))$, the function $t\mapsto\mathcal V_R(u(t))$ is continuous on $I$. The
fundamental theorem of calculus applied to the integrated identity therefore gives
\[
\mathcal M_R(u(\cdot))\in C^1(I),
\quad
\frac{d}{dt}\mathcal M_R(u(t))
=
\mathcal V_R(u(t)).
\]
Hence the differential virial identity remains valid for the nonsmooth exponential weight $\phi_R$. The argument
does not require differentiating the limiting virial observable directly or interchanging the limit
$\delta\downarrow0$ with a time derivative. Instead, we first establish the integrated identity by regularization
and then recover its differential form from the continuity of the limiting virial integrand.

\section{The fractional defect and radial weighted interpolation}
\label{sec:defect}

The fractional virial identity controls a weighted resolvent gradient rather than a localized $\dot H^s$ norm. The purpose of this section is to transfer that control to dyadic fractional Sobolev square functions and then use radial Sobolev decay to estimate the nonlinear localization error.

\subsection{Dyadic decomposition and defect quantities}

Let
\[
\rho_j:=2^jR, \quad A_j:=\{x:\rho_j\le|x|<2\rho_j\}, \quad j\in\mathbb Z.
\]
Choose smooth radial functions $\zeta_j,\eta_j$ as in Appendix~\ref{app:cutoffs} such that
\begin{equation}
\sum_j\zeta_j(x)^2=1
\quad (x\ne0),
\label{eq:quadratic-partition}
\end{equation}
\[
\supp\zeta_j\subset A_j^*:=\{2^{-1}\rho_j<|x|<4\rho_j\},
\]
\[
\eta_j\equiv1\ \text{on }\supp\zeta_j,
\quad
\supp\eta_j\subset A_j^{**}:=\{2^{-2}\rho_j<|x|<8\rho_j\},
\]
and, for $\chi_j\in\{\zeta_j,\eta_j\}$,
\begin{equation}
|\nabla^k\chi_j|\le C_k\rho_j^{-k},
\quad k=1,2.
\label{eq:cutoff-derivatives}
\end{equation}
The families are uniformly locally finite: there exists $C>0$, independent of $R$, such that every $x\ne0$ belongs to the support of at most $C$ functions $\zeta_j$, and likewise for the family $\{\eta_j\}$. Set
\[
d_j:=\min\left\{\left(\frac{\rho_j}{R}\right)^{2s},1\right\}.
\]
On every fixed enlargement of $A_j^*$ or $A_j^{**}$,
\begin{equation}
d_j\simeq\omega_R(x).
\label{eq:dj-omega}
\end{equation}
For an admissible annular family $\chi=\{\chi_j\}$, define
\[
 \mathcal D_R^\chi(u):=\sum_jd_j\|(-\Delta)^{s/2}(\chi_ju)\|_{L^2}^2,
\]
and set
\[
 \mathcal D_R(u):= \mathcal D_R^\zeta(u)+ \mathcal D_R^\eta(u).
\]
Motivated by the weighted kinetic term arising in the fractional virial identity, we introduce the Balakrishnan defect
\[
 \mathcal Q_R(u):=\int_0^\infty\mu^s\int_{\R^N}\omega_R(x)|\nabla u_\mu(x)|^2 dx d\mu.
\]

\subsection{An off-diagonal resolvent estimate}

We begin with an off-diagonal estimate for the resolvent between dyadic annuli.

\begin{lemma}[Dyadic off-diagonal resolvent estimate]
\label{lem:off-diagonal}
Let $\widetilde A_j$ be a fixed enlargement of a dyadic annulus at radius
$\rho_j=2^jR$. Then there exist constants $C,\kappa>0$, depending only on
$N$ and on the chosen enlargement, such that for every $j\in\mathbb Z$,
every $\mu\ge\rho_j^{-2}$, and every $u\in L^2(\R^N)$,
\begin{equation}
\|\1_{\widetilde A_j}(-\Delta+\mu)^{-1}u\|_{L^2}^2
\le
C\mu^{-2}\sum_{k\in\mathbb Z}2^{-\kappa|j-k|}\|\zeta_ku\|_{L^2}^2.
\label{eq:off-diagonal}
\end{equation}
\end{lemma}

\begin{proof}
Write
\[
G_\mu:=(-\Delta+\mu)^{-1}.
\]
Since $\sum_k\zeta_k^2=1$ on $\R^N\setminus\{0\}$, we have in $L^2(\R^N)$
\[
u=\sum_{k\in\mathbb Z}\zeta_k^2u.
\]
Hence
\[
\1_{\widetilde A_j}G_\mu u
=
\sum_k\1_{\widetilde A_j}G_\mu(\zeta_k^2u).
\]
Because $\zeta_k^2u$ is supported in $\supp\zeta_k$,
\[
G_\mu(\zeta_k^2u)
=
G_\mu\1_{\supp\zeta_k}(\zeta_k^2u).
\]
It is therefore enough to prove a pairwise bound
\begin{equation}
\|\1_{\widetilde A_j}G_\mu\1_{\supp\zeta_k}\|_{L^2\to L^2}
\le
C\mu^{-1}2^{-\kappa_0|j-k|}
\label{eq:pairwise-resolvent}
\end{equation}
for some $\kappa_0>0$. In fact, by \eqref{eq:pairwise-resolvent}, 
\[
\begin{aligned}
\|\1_{\widetilde A_j}G_\mu u\|_{L^2}
&\le
\sum_k
\|\1_{\widetilde A_j}G_\mu\1_{\supp\zeta_k}(\zeta_k^2u)\|_{L^2}\\
&\lesssim
\mu^{-1}\sum_k
2^{-\kappa_0|j-k|}\|\zeta_k^2u\|_{L^2}\\
&\le
C\mu^{-1}\sum_k
2^{-\kappa_0|j-k|}\|\zeta_ku\|_{L^2}.
\end{aligned}
\]
Since $\sum_k2^{-\kappa_0|j-k|}$ is uniformly bounded, Cauchy--Schwarz inequality gives
\[
\|\1_{\widetilde A_j}G_\mu u\|_{L^2}^2
\lesssim
\mu^{-2}
\sum_k2^{-\kappa_0|j-k|}\|\zeta_ku\|_{L^2}^2.
\]
This is \eqref{eq:off-diagonal}.

To prove \eqref{eq:pairwise-resolvent}, we use the heat-semigroup representation
\begin{equation}
G_\mu=\int_0^\infty e^{-\mu t}e^{t\Delta} dt,
\label{eq:heat-resolvent}
\end{equation}
where the integral converges strongly on $L^2(\R^N)$. We also use the Davies--Gaffney estimate
\begin{equation}
\|\1_Ee^{t\Delta}\1_F\|_{L^2\to L^2}
\le
\exp\left(-\frac{\dist(E,F)^2}{4t}\right),
\quad t>0,
\label{eq:Davies-Gaffney}
\end{equation}
for measurable sets $E,F\subset\R^N$; see, for example, \cite[Theorem 2]{Davies}. In the Euclidean setting it follows directly from the Gaussian heat kernel
\begin{equation} \label{eq:heat-kernel}
K_t(x,y) = (4\pi t)^{-N/2} \exp\left(-\frac{|x-y|^2}{4t}\right).
\end{equation}

We distinguish three cases.

\medskip
\noindent\emph{Case 1: $k\ge j+C_0$.}
Choose $C_0$ sufficiently large, depending only on the fixed enlargement $\widetilde A_j$. Since
\[
|x|\simeq\rho_j \quad\text{on }\widetilde A_j,
\quad
|x|\simeq\rho_k \quad\text{on }\supp\zeta_k,
\]
we have
\[
\dist(\widetilde A_j,\supp\zeta_k)\simeq\rho_k.
\]
By \eqref{eq:heat-resolvent} and \eqref{eq:Davies-Gaffney},
\[
\|\1_{\widetilde A_j}G_\mu\1_{\supp\zeta_k}\|_{L^2 \to L^2}
\lesssim
\int_0^\infty e^{-\mu t-c\rho_k^2/t} dt.
\]
With the change of variables $\tau=\mu t$,
\[
\int_0^\infty e^{-\mu t-c\rho_k^2/t} dt = \mu^{-1}\int_0^\infty e^{-\tau-c\mu\rho_k^2/\tau} d\tau.
\]
Set $\lambda:=\sqrt\mu \rho_k$. Since $\mu\ge\rho_j^{-2}$,
\[
\lambda = \sqrt\mu \rho_j\frac{\rho_k}{\rho_j} \ge 2^{k-j}.
\]
For every $M>0$,
\begin{equation}
\int_0^\infty e^{-\tau-c\lambda^2/\tau} d\tau \le C_M(1+\lambda)^{-M}.
\label{eq:outer-decay}
\end{equation}
Indeed, for $\lambda\ge1$ split the integral into $(0,\lambda)$ and $(\lambda,\infty)$. On the first interval $\lambda^2/\tau\ge\lambda$, while on the second $\tau\ge\lambda$; hence the integral is bounded by $C(1+\lambda)e^{-c'\lambda}$, which implies \eqref{eq:outer-decay}. Thus
\[
\|\1_{\widetilde A_j}G_\mu\1_{\supp\zeta_k}\|_{L^2 \to L^2} \le C_M\mu^{-1}2^{-M(k-j)}.
\]

\medskip
\noindent\emph{Case 2: $k\le j-C_0$.}
Now
\[
\dist(\widetilde A_j,\supp\zeta_k)\simeq\rho_j.
\]
At the lowest admissible frequency $\mu\simeq\rho_j^{-2}$, \eqref{eq:Davies-Gaffney} alone gives no decay in $j-k$, so we also
use the volume difference between the two annuli. The Gaussian kernel \eqref{eq:heat-kernel} gives the Hilbert--Schmidt estimate
\[
\|\1_Ee^{t\Delta}\1_F\|_{L^2\to L^2}
\le \|\1_Ee^{t\Delta}\1_F\|_{\mathrm{HS}}
\lesssim t^{-N/2}|E|^{1/2}|F|^{1/2} e^{-c\dist(E,F)^2/t}.
\]
Indeed, using $|x-y|\geq \dist(E,F)$ for all $x\in E, y\in F$, \eqref{eq:heat-kernel} gives
\[
\|\1_Ee^{t\Delta}\1_F\|_{\mathrm{HS}}^2
= \int_E\int_F|K_t(x,y)|^2 dy dx \lesssim t^{-N}|E||F|e^{-c\dist(E,F)^2/t}.
\]
Taking $E=\widetilde A_j$ and $F=\supp\zeta_k$, and using $|E|\simeq\rho_j^N$ and $|F|\lesssim\rho_k^N$, we obtain
\[
\|\1_{\widetilde A_j}e^{t\Delta}\1_{\supp\zeta_k}\|_{L^2 \to L^2}
\lesssim
t^{-N/2}(\rho_j\rho_k)^{N/2}e^{-c\rho_j^2/t}.
\]
Consequently,
\begin{equation}
\|\1_{\widetilde A_j}G_\mu\1_{\supp\zeta_k}\|_{L^2 \to L^2}
\lesssim
(\rho_j\rho_k)^{N/2} \int_0^\infty t^{-N/2}e^{-\mu t-c\rho_j^2/t} dt.
\label{eq:inner-integral}
\end{equation}
Set $t=\rho_j^2\tau$ and $a:=\mu\rho_j^2\ge1$. The integral in \eqref{eq:inner-integral} equals
\[
\rho_j^{2-N} \int_0^\infty \tau^{-N/2}e^{-a\tau-c/\tau} d\tau.
\]

We claim that
\begin{equation}
\label{eq:inner-elementary}
\int_0^\infty \tau^{-N/2}e^{-a\tau-c/\tau} d\tau \lesssim a^{-1},
\quad a\ge1.
\end{equation}
To estimate this integral uniformly in $a\ge1$, split its exponent as
\[
a\tau+\frac c\tau
=\frac12\left(\tau+\frac c\tau\right)
+\left(a-\frac12\right)\tau+\frac{c}{2\tau}.
\]
By the arithmetic--geometric mean inequality,
\[
\left(a-\frac12\right)\tau+\frac{c}{2\tau}
\ge\sqrt{c(2a-1)}\ge\sqrt{ac},\quad a\ge1.
\]
Consequently,
\[
\begin{aligned}
\int_0^\infty\tau^{-N/2}e^{-a\tau-c/\tau} d\tau &\le e^{-\sqrt{ac}}\int_0^\infty
\tau^{-N/2}e^{-(\tau+c/\tau)/2} d\tau\\
&\le C_{N,c}e^{-\sqrt{ac}} \\
&\le C_{N,c}a^{-1}.
\end{aligned}
\]
Here the last integral converges at $\tau=0$ by the factor $e^{-c/(2\tau)}$ and at $\tau=\infty$ by $e^{-\tau/2}$; moreover, $\sup_{a\ge1}a e^{-\sqrt{ac}}<\infty$. This proves \eqref{eq:inner-elementary}.

Substitution in \eqref{eq:inner-integral} yields
\[
\|\1_{\widetilde A_j}G_\mu\1_{\supp\zeta_k}\|_{L^2 \to L^2} \lesssim \mu^{-1} \left(\frac{\rho_k}{\rho_j}\right)^{N/2} = \mu^{-1}2^{-\frac N2(j-k)}.
\]

\medskip
\noindent\emph{Case 3: $|j-k|\le C_0$.}
No off-diagonal gain is needed. By the spectral theorem,
\[
\|G_\mu\|_{L^2 \to L^2} =\sup_{\xi\in\R^N}\frac1{|\xi|^2+\mu} \le\mu^{-1}.
\]
Combining the three cases, and choosing any $0<\kappa_0<N/2$, gives \eqref{eq:pairwise-resolvent}.

\end{proof}

\subsection{Balakrishnan-to-dyadic defect transfer}

\begin{proposition}[Dyadic Balakrishnan localization]
\label{prop:defect-transfer}
Let $N\ge1$ and $0<s<1$. Let $\chi=\{\chi_j\}_{j\in\mathbb Z}$ denote either of the admissible cutoff families $\{\zeta_j\}$ or $\{\eta_j\}$. Then there exist $C, c>0$ such that for every $u\in H^s(\R^N)$ and every $R>0$,
\begin{equation}
 \mathcal D_R^\chi(u) \le C  \mathcal Q_R(u)+CR^{-2s}M(u).
\label{eq:defect-transfer-family}
\end{equation}
Consequently,
\begin{equation}
 \mathcal D_R(u)\le C  \mathcal Q_R(u)+CR^{-2s}M(u),
\label{eq:defect-transfer}
\end{equation}
and therefore
\begin{equation}
 \mathcal Q_R(u)\ge c \mathcal D_R(u)-CR^{-2s}M(u).
\label{eq:Q-lower}
\end{equation}
\end{proposition}

\begin{proof}
Fix one of the two cutoff families and set
\[
f_j:=\chi_j u, \quad G_\mu:=(-\Delta+\mu)^{-1}.
\]
Let $\widetilde A_j$ be a fixed enlargement of $\supp\chi_j$ containing the supports of $\chi_j$, $\nabla\chi_j$, and $\Delta\chi_j$.
The family $\{\widetilde A_j\}$ is uniformly locally finite.

\medskip
\noindent\emph{Step 1: Low resolvent frequencies.}
Since multiplication by $\chi_j$ is bounded on $H^s(\R^N)$, we have $f_j\in H^s(\R^N)$. By \eqref{eq:balakrishnan-energy},
\begin{equation}
s\|(-\Delta)^{s/2}f_j\|_{L^2}^2=c_s^2\int_0^\infty \mu^s\|\nabla G_\mu f_j\|_{L^2}^2 d\mu.
\label{eq:localized-balakrishnan}
\end{equation}
We split the integral at
\[
\mu=\rho_j^{-2}.
\]
By Plancherel,
\[
\|\nabla G_\mu\|_{L^2 \to L^2}=\sup_{\xi\in\R^N}\frac{|\xi|}{|\xi|^2+\mu} \lesssim\mu^{-1/2}.
\]
Therefore
\[
\|\nabla G_\mu f_j\|_{L^2}^2 \lesssim \mu^{-1}\|f_j\|_{L^2}^2,
\]
and, since $s>0$,
\begin{equation}
\int_0^{\rho_j^{-2}} \mu^s\|\nabla G_\mu f_j\|_{L^2}^2 d\mu \lesssim \rho_j^{-2s}\|\chi_j u\|_{L^2}^2.
\label{eq:low-resolvent}
\end{equation}

\medskip
\noindent\emph{Step 2: High resolvent frequencies.}
Assume now that $\mu\ge\rho_j^{-2}$ and put
\[
v_\mu:=G_\mu u.
\]
Since $(-\Delta+\mu)v_\mu=u$,
\[
(-\Delta+\mu)(\chi_jv_\mu)=\chi_j u-[\Delta,\chi_j]v_\mu,
\]
where
\[
[\Delta,\chi_j]v_\mu = 2\nabla\chi_j\cdot\nabla v_\mu+(\Delta\chi_j)v_\mu.
\]
Applying $G_\mu$ gives, up to the harmless sign convention for the commutator,
\begin{equation}
G_\mu(\chi_j u)=\chi_jG_\mu u+G_\mu[\Delta,\chi_j]G_\mu u.
\label{eq:commutator-resolvent}
\end{equation}
Differentiating \eqref{eq:commutator-resolvent},
\[
\nabla G_\mu f_j=\nabla(\chi_jv_\mu)+\nabla G_\mu
\bigl(2\nabla\chi_j\cdot\nabla v_\mu+(\Delta\chi_j)v_\mu\bigr).
\]

For the first term,
\[
\nabla(\chi_jv_\mu)=\chi_j\nabla v_\mu+(\nabla\chi_j)v_\mu,
\]
so \eqref{eq:cutoff-derivatives} yields
\begin{equation}
\|\nabla(\chi_jv_\mu)\|_{L^2}^2 \lesssim \|\1_{\widetilde A_j}\nabla v_\mu\|_{L^2}^2
+\rho_j^{-2}\|\1_{\widetilde A_j}v_\mu\|_{L^2}^2.
\label{eq:main-local-resolvent}
\end{equation}
For the commutator term, the bound
$\|\nabla G_\mu\|_{L^2 \to L^2}\lesssim\mu^{-1/2}$ gives
\[
\begin{aligned}
\|\nabla G_\mu[\Delta,\chi_j]v_\mu\|_{L^2}^2\lesssim \mu^{-1}\rho_j^{-2} \|\1_{\widetilde A_j}\nabla v_\mu\|_{L^2}^2+\mu^{-1}\rho_j^{-4} \|\1_{\widetilde A_j}v_\mu\|_{L^2}^2.
\end{aligned}
\]
Combining this with \eqref{eq:main-local-resolvent},
\begin{align}
\|\nabla G_\mu f_j\|_{L^2}^2 \lesssim{}& \|\1_{\widetilde A_j}\nabla G_\mu u\|_{L^2}^2
+ \rho_j^{-2}\|\1_{\widetilde A_j}G_\mu u\|_{L^2}^2 \notag\\
&+\mu^{-1}\rho_j^{-2} \|\1_{\widetilde A_j}\nabla G_\mu u\|_{L^2}^2
+\mu^{-1}\rho_j^{-4}\|\1_{\widetilde A_j}G_\mu u\|_{L^2}^2.
\label{eq:high-resolvent-pre}
\end{align}
Since $\mu\ge\rho_j^{-2}$, $\mu^{-1}\rho_j^{-2}\le1$, so the third term is absorbed into the first.

\medskip
\noindent\emph{Step 3: Off-diagonal resolvent control.}
Lemma \ref{lem:off-diagonal} gives, for $\mu\ge\rho_j^{-2}$,
\[
\|\1_{\widetilde A_j}G_\mu u\|_{L^2}^2 \lesssim \mu^{-2}\sum_k2^{-\kappa|j-k|}\|\zeta_ku\|_{L^2}^2.
\]
Substitution into \eqref{eq:high-resolvent-pre} yields
\begin{align}
\|\nabla G_\mu f_j\|_{L^2}^2 \lesssim{}& \|\1_{\widetilde A_j}\nabla G_\mu u\|_{L^2}^2 \notag\\
&+\rho_j^{-2}\mu^{-2} \sum_k2^{-\kappa|j-k|}\|\zeta_ku\|_{L^2}^2 \notag\\
&+\rho_j^{-4}\mu^{-3} \sum_k2^{-\kappa|j-k|}\|\zeta_ku\|_{L^2}^2.
\notag
\end{align}
Multiplying by $\mu^s$ and integrating over $[\rho_j^{-2},\infty)$, we use $s<1$ and
\[
\rho_j^{-2}\int_{\rho_j^{-2}}^\infty\mu^{s-2} d\mu=C_s\rho_j^{-2s},
\]
as well as
\[
\rho_j^{-4}\int_{\rho_j^{-2}}^\infty\mu^{s-3} d\mu=C_s\rho_j^{-2s},
\]
to obtain
\begin{align}
\int_{\rho_j^{-2}}^\infty \mu^s\|\nabla G_\mu f_j\|_{L^2}^2 d\mu \lesssim 
\int_{\rho_j^{-2}}^\infty \mu^s\|\1_{\widetilde A_j}\nabla G_\mu u\|_{L^2}^2 d\mu + \rho_j^{-2s}\sum_k 2^{-\kappa|j-k|}\|\zeta_ku\|_{L^2}^2.
\label{eq:high-resolvent}
\end{align}
Combining \eqref{eq:localized-balakrishnan}, \eqref{eq:low-resolvent}, and \eqref{eq:high-resolvent}, and enlarging the principal $\mu$-integral to $(0,\infty)$, gives
\begin{align}
\|(-\Delta)^{s/2}(\chi_j u)\|_{L^2}^2\lesssim{}&
\int_0^\infty \mu^s\|\1_{\widetilde A_j}\nabla G_\mu u\|_{L^2}^2 d\mu \notag\\
&+\rho_j^{-2s}\|\chi_j u\|_{L^2}^2
+\rho_j^{-2s}\sum_k 2^{-\kappa|j-k|}\|\zeta_ku\|_{L^2}^2.
\label{eq:localized-Hs}
\end{align}

\medskip
\noindent\emph{Step 4: Weighted summation.}
Multiply \eqref{eq:localized-Hs} by $d_j$ and sum over $j$. On $\widetilde A_j$, \eqref{eq:dj-omega} gives
$d_j\lesssim\omega_R(x)$. Since the enlarged annuli are uniformly locally finite,
\[
\sum_jd_j\1_{\widetilde A_j}(x) \lesssim\omega_R(x).
\]
Tonelli's theorem therefore yields
\[
\begin{aligned}
\sum_jd_j\int_0^\infty \mu^s\|\1_{\widetilde A_j}\nabla G_\mu u\|_{L^2}^2 d\mu \lesssim \int_0^\infty\mu^s \int_{\R^N}\omega_R(x)|\nabla G_\mu u(x)|^2 dx d\mu.
\end{aligned}
\]
Since $u_\mu=c_sG_\mu u$, the right-hand side is $c_s^{-2} \mathcal Q_R(u)$ and hence is bounded by $C  \mathcal Q_R(u)$.

Observe that for every $j$,
\begin{equation}
d_j\rho_j^{-2s}\le R^{-2s}.
\label{eq:dj-scale}
\end{equation}
Indeed, equality holds when $\rho_j\le R$, whereas for $\rho_j\ge R$ we have $d_j=1$ and $\rho_j^{-2s}\le R^{-2s}$. Therefore
\[
\sum_jd_j\rho_j^{-2s}\|\chi_j u\|_{L^2}^2
\lesssim R^{-2s}\|u\|_{L^2}^2= R^{-2s}M(u),
\]
where we used the uniform local finiteness of the family $\{\chi_j\}$. For $\chi=\zeta$, the quadratic partition in fact gives
$\sum_j\|\zeta_ju\|_{L^2}^2=\|u\|_{L^2}^2$ exactly.

Using \eqref{eq:dj-scale} and Tonelli's theorem,
\[
\begin{aligned}
&\sum_jd_j\rho_j^{-2s} \sum_k2^{-\kappa|j-k|}\|\zeta_ku\|_{L^2}^2\\
&\quad\le R^{-2s} \sum_k\|\zeta_ku\|_{L^2}^2 \sum_j2^{-\kappa|j-k|}
\lesssim R^{-2s}M(u),
\end{aligned}
\]
because $\sum_{\ell\in\mathbb Z}2^{-\kappa|\ell|}<\infty$ and $\sum_k\|\zeta_ku\|_{L^2}^2=M(u)$.

Combining the above inequalities gives
\[
 \mathcal D_R^\chi(u) \lesssim  \mathcal Q_R(u)+R^{-2s}M(u),
\]
which is \eqref{eq:defect-transfer-family}. Applying the estimate to $\chi=\zeta$ and $\chi=\eta$ and adjusting the constant gives
\eqref{eq:defect-transfer}. Rearranging then yields \eqref{eq:Q-lower}.
\end{proof}

\subsection{Radial weighted interpolation}

\begin{proposition}[Radial weighted nonlinear interpolation]
\label{prop:weighted-interpolation}
Assume \eqref{eq:range}. Let $\varepsilon>0$. There exists $C_\varepsilon>0$ such that for every radial $u\in H^s(\R^N)$ and every $R>0$,
\begin{equation}
\int_{\R^N}\omega_R|u|^q dx \le \varepsilon  \mathcal D_R(u)+C_\varepsilon R^{-2s}M(u)^{N/(N-2s)}.
\label{eq:weighted-nonlinear-young}
\end{equation}
\end{proposition}

\begin{proof}
Set
\[
\theta:=\frac{2s}{N}, \quad \beta:=2s(1-\theta)=\frac{2s(N-2s)}{N}.
\]
By \eqref{eq:range},
\[
0<\theta<1.
\]
It suffices to prove
\begin{equation}
\int_{\R^N}\omega_R(x)|u(x)|^q dx \le CR^{-\beta}M(u)  \mathcal D_R(u)^\theta.
\label{eq:weighted-nonlinear}
\end{equation}
In fact, if \eqref{eq:weighted-nonlinear} holds, then Young's inequality with exponents $1/\theta$ and $1/(1-\theta)$ gives
\[
R^{-\beta}M(u) \mathcal D_R(u)^\theta
\le \varepsilon  \mathcal D_R(u) +C_\varepsilon R^{-\beta/(1-\theta)}M(u)^{1/(1-\theta)},
\]
which is \eqref{eq:weighted-nonlinear-young} because
\[
\frac{\beta}{1-\theta}=2s, \quad \frac1{1-\theta}=\frac{N}{N-2s}.
\]

We now prove \eqref{eq:weighted-nonlinear}. Put
\[
m_j:=\|\zeta_ju\|_{L^2}^2, \quad K_j:=\|(-\Delta)^{s/2}(\eta_ju)\|_{L^2}^2, \quad D_j:=d_jK_j.
\]
Since $\eta_j\equiv1$ on $\supp\zeta_j$ and $\eta_ju$ is radial, \eqref{eq:radial-sobolev} gives for $x\in\supp\zeta_j$
\begin{equation}
|u(x)|=|\eta_j(x)u(x)| \lesssim \rho_j^{-(N/2-s)}K_j^{1/2}.
\label{eq:annular-radial}
\end{equation}
Using \eqref{eq:quadratic-partition}, $q-2=4s/N$, and $\omega_R\simeq d_j$ on $\supp\zeta_j$,
\begin{align}
\int\omega_R|u|^q dx
&=\sum_j\int\zeta_j^2\omega_R|u|^q dx \notag\\
&\lesssim\sum_jd_j\rho_j^{-\beta}K_j^\theta m_j
=\sum_j\rho_j^{-\beta}d_j^{1-\theta}D_j^\theta m_j.
\label{eq:annular-nonlinear-sum}
\end{align}
The choice $\beta=2s(1-\theta)$ gives
\begin{equation}
\rho_j^{-\beta}d_j^{1-\theta}\lesssim R^{-\beta}
\label{eq:scale-match}
\end{equation}
uniformly in $j$. Indeed, if $\rho_j\le R$, then the left-hand side is comparable to $R^{-\beta}$, while for $\rho_j\ge R$ it is bounded by $R^{-\beta}$.

By H\"older's inequality in $j$,
\begin{align*}
\sum_jD_j^\theta m_j
&\le \left(\sum_jD_j\right)^\theta
\left(\sum_jm_j^{1/(1-\theta)}\right)^{1-\theta}.
\end{align*}
Since $1/(1-\theta)>1$, the embedding $\ell^1\subset\ell^{1/(1-\theta)}$ gives
\[
\left(\sum_jm_j^{1/(1-\theta)}\right)^{1-\theta} \le \sum_jm_j = M(u),
\]
where the last identity follows from \eqref{eq:quadratic-partition}. Moreover, $\sum_jD_j= \mathcal D_R^\eta(u)\le  \mathcal D_R(u)$. Hence
\[
\sum_jD_j^\theta m_j \le  \mathcal D_R(u)^\theta M(u).
\]
Together with \eqref{eq:annular-nonlinear-sum}--\eqref{eq:scale-match}, this proves \eqref{eq:weighted-nonlinear}.
\end{proof}

\begin{remark}[Where radiality enters]
\label{rem:radiality}
The radial assumption is used only in the annular estimate \eqref{eq:annular-radial}. The geometry of the multiplier, the Balakrishnan defect transfer, the virial--defect estimate below, and the final Osgood argument are otherwise nonradial. Thus a nonradial substitute for \eqref{eq:weighted-nonlinear-young} with an arbitrarily small coefficient in front of $ \mathcal D_R$ would extend the finite-time theorem beyond radial data.
\end{remark}

\section{Localized fractional virial inequality}
\label{sec:virial}

For the exponential weight introduced in Section \ref{sec:weight}, define
\[
 \mathcal M_R(t):=2\operatorname{Im}\int_{\R^N}\overline{u(t,x)} \nabla\phi_R(x)\cdot\nabla u(t,x) dx.
\]
Since $s>1/2$, this quadratic form extends continuously to $H^s(\R^N)$; see Lemma \ref{lem:first-order-form} below.

We first estimate the lower-order term generated by $\Delta^2\phi_R$.

\begin{lemma}[Lower-order Balakrishnan term]
\label{lem:lower-balakrishnan}
Let $N\ge2$ and $0<s<1$. There exists $C>0$ such that for every $u\in L^2(\R^N)$ and every $R>0$,
\begin{equation}
\left|
\int_0^\infty\mu^s\int_{\R^N}(\Delta^2\phi_R)|u_\mu|^2 dx d\mu
\right|
\le CR^{-2s} M(u).
\label{eq:lower-balakrishnan}
\end{equation}
\end{lemma}

\begin{proof}
By \eqref{eq:bilap-weight} and the fractional Hardy inequality \eqref{eq:hardy} with $\alpha=1-s$,
\begin{align*}
\left|
\int_0^\infty\mu^s\int(\Delta^2\phi_R)|u_\mu|^2 dx d\mu
\right|
&\lesssim R^{-2s} \int_0^\infty\mu^s\||x|^{-(1-s)}u_\mu\|_{L^2}^2 d\mu\\
&\lesssim R^{-2s} \int_0^\infty\mu^s \|(-\Delta)^{(1-s)/2}u_\mu\|_{L^2}^2 d\mu.
\end{align*}
By Plancherel and \eqref{eq:umu}, the last integral is $C_s\|u\|_{L^2}^2$ (see \eqref{eq:iden-L2}). This proves \eqref{eq:lower-balakrishnan}.
\end{proof}

\begin{proposition}[Coercive localized virial inequality]
\label{prop:coercive-virial}
Assume \eqref{eq:range}. Let
\[
u\in C((-T_-,T_+);H^{2s}_{\rm rad}(\R^N))
\]
be a radial strong solution of \eqref{eq:fnls}. There exist constants $c_{N,s}, C_{N,s}>0$ such that for every $R>0$ and every $t\in (-T_-,T_+)$,
\begin{equation}
\mathcal M_R'(t) \le 8sE(u_0)-c_{N,s}  \mathcal D_R(u(t)) +C_{N,s} R^{-2s}\left[m+m^{N/(N-2s)}\right],
\label{eq:coercive-virial}
\end{equation}
where $m=M(u_0)$.
\end{proposition}

\begin{proof}
Apply \eqref{eq:BHL-virial} with $\phi=\phi_R$. The principal kinetic term is
\begin{align*}
4\int_0^\infty\mu^s\int\partial_k\overline{u_\mu}(\partial_{kl}\phi_R)\partial_lu_\mu dx d\mu
={}&4\int_0^\infty\mu^s\|\nabla u_\mu\|_{L^2}^2 d\mu\\
&-4\int_0^\infty\mu^s\int \nabla\overline{u_\mu}\cdot(I-D^2\phi_R)\nabla u_\mu dx d\mu.
\end{align*}
By \eqref{eq:balakrishnan-energy}, \eqref{eq:hessian-defect}, and Proposition \ref{prop:defect-transfer},
\begin{equation}
\text{kinetic term}
\le 4s\|(-\Delta)^{s/2}u\|_{L^2}^2 -c \mathcal D_R(u)+CR^{-2s} M(u).
\label{eq:kinetic-main}
\end{equation}
The lower-order term is controlled by Lemma \ref{lem:lower-balakrishnan}.

For the nonlinearity, we use $p=\frac{4s}{N}$ and $q=p+2$ (see \eqref{eq:p-q}). Then
\begin{align*}
-\frac{2p}{p+2}\int(\Delta\phi_R)|u|^{p+2} dx =-\frac{8s}{q}\int|u|^q dx+\frac{2p}{q}\int(N-\Delta\phi_R)|u|^q dx.
\end{align*}
Consequently, the leading quadratic contribution in \eqref{eq:kinetic-main} and the leading nonlinear term combine exactly into
\[
4s\|(-\Delta)^{s/2}u\|_{L^2}^2-\frac{8s}{q}\|u\|_{L^q}^q =8sE(u)=8sE(u_0).
\]
It remains to estimate the positive localization error. By \eqref{eq:lap-defect} and Proposition \ref{prop:weighted-interpolation}, for every $\varepsilon>0$,
\[
\int(N-\Delta\phi_R)|u|^q dx \lesssim \varepsilon  \mathcal D_R(u)
+C_\varepsilon R^{-2s} M(u)^{N/(N-2s)}.
\]
Choose $\varepsilon$ sufficiently small to absorb this term into the negative defect contribution in \eqref{eq:kinetic-main}, and use mass conservation. This proves \eqref{eq:coercive-virial}.
\end{proof}

\begin{proposition}[Integrated virial inequality for energy-space solutions]
\label{prop:energy-virial}
Assume \eqref{eq:expo-energy} and let $u\in C((-T_-,T_+);H^s_{\rm rad}(\R^N))$ be its maximal solution with $u(0)=u_0$. For every $R>0$ and $-T_-<t_1\le t_2<T_+$, one has
\begin{align}
\mathcal M_R(u(t_2))-\mathcal M_R(u(t_1))
\le{}&8s(t_2-t_1)E(u_0)
-c_{N,s}\int_{t_1}^{t_2}\mathcal D_R(u(t)) dt\notag\\
&+C_{N,s}(t_2-t_1)R^{-2s}
\left[M(u_0)+M(u_0)^{\frac{N}{N-2s}}\right].
\label{eq:energy-virial-integrated}
\end{align}
The integral of $\mathcal D_R$ is finite on compact subintervals.
\end{proposition}

\begin{proof}
The proof is given in Appendix~\ref{app:energy-approx}, where we establish the integrated virial identity for radial $H^s(\R^N)$ solutions and apply the estimates underlying Proposition~\ref{prop:coercive-virial} directly to its integrand. This allows us to obtain the desired inequality without assuming $H^{2s}$ regularity. The appendix also justifies the finiteness of the defect integral on every
compact subinterval of the solution's lifespan.
\end{proof}

\section{The virial--defect estimate}
\label{sec:virial-defect}

We now show that the same dyadic defect appearing with a negative sign in Proposition \ref{prop:coercive-virial} also controls the localized virial observable.

\begin{proposition}[Virial--defect estimate]
\label{prop:virial-defect}
Let $N\ge 1$ and $\frac12 < s<1$. There exist $C_0, C_1>0$ depending only on $N$ and $s$ such that for every $u\in H^s(\R^N)$ and every $R>0$,
\begin{equation}
| \mathcal M_R(u)| \le C_0 M(u)+ C_1R M(u)^{1-\frac1{2s}} \mathcal D_R(u)^{\frac1{2s}}.
\label{eq:virial-defect}
\end{equation}
\end{proposition}

To prove Proposition \ref{prop:virial-defect}, we use the following refinement of \cite[Lemma~A.1]{BHL} whose proof is given in Appendix \ref{app:first-order-form}.

\begin{lemma}[First-order quadratic form]
\label{lem:first-order-form}
Let $a:\R^N\to\R^N$ be a real-valued bounded Lipschitz vector field. Then there exists $C_N>0$ such that for every $f\in H^{1/2}(\R^N)$,
\begin{equation}
\left| 2\operatorname{Im}\int_{\R^N}\overline f a\cdot\nabla f dx \right|
\le C_N\|a\|_{L^\infty}\|(-\Delta)^{1/4}f\|_{L^2}^2
+ C_N\|\nabla a\|_{L^\infty}\|f\|_{L^2}^2.
\label{eq:first-order-form}
\end{equation}
For $f\in H^{1/2}(\R^N)$, the quadratic form on the left-hand side is understood as its continuous extension from Schwartz functions.
\end{lemma}

\begin{remark} \label{rem:density-H1/2}
For a general function $f\in H^{1/2}(\R^N)$, the integral in \eqref{eq:first-order-form} need not be well-defined in the usual sense, since $\nabla f$ may not belong to $L^2$. We therefore define the quadratic form by approximating $f$ with Schwartz functions. The estimate in Lemma~\ref{lem:first-order-form} yields
\[
|Q_a(f)-Q_a(g)|
\lesssim_a
\bigl(\|f\|_{H^{1/2}}+\|g\|_{H^{1/2}}\bigr)
\|f-g\|_{H^{1/2}}
\]
for Schwartz functions $f$ and $g$. Consequently, the quadratic-form values of any sequence converging in $H^{1/2}$ form a Cauchy sequence, and the resulting limit is independent of the approximation. This defines a unique continuous extension of the
quadratic form to $H^{1/2}(\R^N)$.
\end{remark}

\begin{proof}[Proof of Proposition \ref{prop:virial-defect}]
We first prove the estimate for $u\in\mathcal S(\R^N)$ and then pass to $H^s(\R^N)$ by density.

Set
\[
f_j:=\zeta_j u, \quad m_j:=\|f_j\|_{L^2}^2, \quad K_j:=\|(-\Delta)^{s/2}f_j\|_{L^2}^2, \quad D_j:=d_jK_j.
\]

\medskip
\noindent\emph{Step 1: Dyadic decomposition of the virial form.}
Using $\sum_j\zeta_j^2=1$,
\[
\mathcal M_R(u) =\sum_j 2\operatorname{Im}\int \zeta_j^2\overline u \nabla\phi_R\cdot\nabla u dx.
\]
Since
\[
\nabla f_j=\zeta_j\nabla u+u\nabla\zeta_j,
\]
we have
\[
\overline{f_j}\nabla f_j =\zeta_j^2\overline u \nabla u+\zeta_j|u|^2\nabla\zeta_j.
\]
The second term is real-valued, so it disappears after taking the imaginary part. Therefore
\[
\mathcal M_R(u)= \sum_j 2\operatorname{Im}\int_{\R^N} \overline{f_j} \nabla\phi_R\cdot\nabla f_j dx.
\]

\medskip
\noindent\emph{Step 2: Localization of the vector field.} Choose fixed constants $0<c<C<\infty$ and set
\[
\widetilde A_j:=\{x\in\R^N:c\rho_j<|x|<C\rho_j\},
\]
so that
\[
\supp\zeta_j\Subset \widetilde A_j.
\]
Choose a radial cutoff $\vartheta_j$ satisfying
\[
0\le\vartheta_j\le1,\quad
\vartheta_j\equiv1\ \text{on }\supp\zeta_j,\quad
\supp\vartheta_j\subset\widetilde A_j,
\]
and
\[
|\nabla\vartheta_j|\lesssim\rho_j^{-1}.
\]
Define
\[
a^{(j)}:=\vartheta_j\nabla\phi_R,
\quad
a_j:=\sup_{\widetilde A_j}|\nabla\phi_R|.
\]
Then $a^{(j)}=\nabla\phi_R$ on $\supp f_j$ and
\[
\|a^{(j)}\|_{L^\infty}\le a_j.
\]
Moreover,
\[
\nabla a^{(j)} = \nabla\vartheta_j\otimes\nabla\phi_R +\vartheta_jD^2\phi_R.
\]
Since $\|D^2\phi_R\|_{L^\infty}\lesssim1$ and $|\nabla\phi_R(x)|\le|x|\lesssim\rho_j$ on $\supp\vartheta_j\subset\widetilde A_j$, we obtain
\[
\|\nabla a^{(j)}\|_{L^\infty}\lesssim1
\]
uniformly in $j$ and $R$.

Applying Lemma \ref{lem:first-order-form} to $f_j$ and $a^{(j)}$,
\begin{equation}
\left| 2\operatorname{Im}\int \overline{f_j} \nabla\phi_R\cdot\nabla f_j dx \right|
\lesssim a_j\|(-\Delta)^{1/4}f_j\|_{L^2}^2+m_j.
\label{eq:local-virial-form}
\end{equation}

\medskip
\noindent\emph{Step 3: Interpolate between $L^2$ and $\dot H^s$.}
Because $s>1/2$, H\"older's inequality on the Fourier side gives
\begin{align}
\|(-\Delta)^{1/4}f_j\|_{L^2}^2
&= \int_{\R^N}|\xi||\widehat f_j(\xi)|^2 d\xi \notag\\
&\le \left(\int|\widehat f_j|^2\right)^{1-\frac1{2s}} \left(\int|\xi|^{2s}|\widehat f_j|^2\right)^{\frac1{2s}} \notag\\
&= m_j^{1-\frac1{2s}}K_j^{\frac1{2s}}.
\label{eq:Hhalf-interp}
\end{align}

\medskip
\noindent\emph{Step 4: Insert the defect weight.} Since $D_j=d_jK_j$,
\[
K_j^{1/(2s)}=d_j^{-1/(2s)}D_j^{1/(2s)}.
\]
Because $\widetilde A_j$ is a fixed enlargement of the dyadic annulus at scale $\rho_j$, we have
\[
\omega_R(x)\simeq d_j \quad\text{for }x\in\widetilde A_j.
\]
Hence, by \eqref{eq:gradient-defect},
\[
|\nabla\phi_R(x)| \lesssim R \omega_R(x)^{1/(2s)} \lesssim R d_j^{1/(2s)}, \quad x\in\widetilde A_j.
\]
Taking the supremum over $\widetilde A_j$ gives
\begin{equation} \label{eq:aj-dj}
a_jd_j^{-1/(2s)}\lesssim R.
\end{equation}
Therefore, combining \eqref{eq:local-virial-form} with the interpolation estimate,
\[
|\mathcal M_R(u)| \lesssim \sum_j m_j + R\sum_j m_j^{1-\frac1{2s}}D_j^{\frac1{2s}}.
\]

\medskip
\noindent\emph{Step 5: Sum over the annuli.}
The quadratic partition gives
\[
\sum_jm_j=\sum_j\|\zeta_ju\|_{L^2}^2=M(u).
\]
Applying H\"older's inequality in $j$ with conjugate exponents $2s/(2s-1)$ and $2s$,
\[
\begin{aligned}
\sum_jm_j^{1-\frac1{2s}}D_j^{\frac1{2s}}
&\le \left(\sum_jm_j\right)^{1-\frac1{2s}} \left(\sum_jD_j\right)^{\frac1{2s}}\\
&= M(u)^{1-\frac1{2s}} \bigl( \mathcal D_R^\zeta(u)\bigr)^{\frac1{2s}}\\
&\le M(u)^{1-\frac1{2s}} \mathcal D_R(u)^{\frac1{2s}}.
\end{aligned}
\]
This proves \eqref{eq:virial-defect} for Schwartz functions.

\medskip
\noindent\emph{Step 6: Extension to $H^s(\R^N)$.}
Let $u_n\in\mathcal S(\R^N)$ with $u_n\to u$ in $H^s(\R^N)$. Since $s>1/2$, Lemma \ref{lem:first-order-form} implies $\mathcal M_R(u_n)\to  \mathcal M_R(u)$, while $M(u_n)\to M(u)$. It remains only to note that $ \mathcal D_R(u_n)\to  \mathcal D_R(u)$. Indeed, Proposition \ref{prop:defect-transfer}, applied to $w=u_n-u$, gives
\[
 \mathcal D_R(w) \lesssim  \mathcal Q_R(w)+R^{-2s}M(w) \lesssim_R \|w\|_{H^s}^2,
\]
because $0\le\omega_R\le1$ and \eqref{eq:balakrishnan-energy} imply $ \mathcal Q_R(w)\lesssim\|(-\Delta)^{s/2}w\|_{L^2}^2$.
Since $ \mathcal D_R^{1/2}$ is the norm of the corresponding weighted $\ell^2(L^2)$ square function,
\[
|\mathcal D_R(u_n)^{1/2}- \mathcal D_R(u)^{1/2}| \le  \mathcal D_R(u_n-u)^{1/2}\longrightarrow0.
\]
Passing to the limit in the Schwartz estimate proves \eqref{eq:virial-defect} for all $u\in H^s(\R^N)$.
\end{proof}

\begin{remark}[Why the order $2s$ is natural]
\label{rem:order}
Suppose instead that the multiplier defect behaved near the origin like
\[
\omega_R(r)\sim \left(\frac rR\right)^\gamma.
\]
There are two competing requirements on $\gamma$.

First, in the virial--defect estimate, the analogue of
\[
|\nabla\phi_R|\omega_R^{-1/(2s)}
\]
behaves near the origin like
\[
r\left(\frac rR\right)^{-\gamma/(2s)}.
\]
Uniform boundedness as $r\downarrow0$ therefore requires $\gamma\le 2s$.

On the other hand, in the dyadic Balakrishnan localization, the low-frequency and commutator mass errors contain the factor
\[
d_j\rho_j^{-2s}, \quad d_j\sim\left(\frac{\rho_j}{R}\right)^\gamma
\quad\text{when }\rho_j\le R.
\]
Hence
\[
d_j\rho_j^{-2s} \sim R^{-\gamma}\rho_j^{\gamma-2s}.
\]
For this quantity to be bounded uniformly in the inner dyadic scales $\rho_j\downarrow0$ by a multiple of $R^{-2s}$, one must have $\gamma\ge 2s$.

Combining both parts requires
\[
\gamma=2s.
\]
This explains why the defect order in the exponential multiplier is matched exactly to the order of the fractional dispersion.
\end{remark}

\section{Proof of finite-time blow-up}
\label{sec:blowup}

\begin{proof}[Proof of Theorem \ref{thm:strong}]
Write
\[
E(u_0)=-e, \quad e>0,
\]
and let $m=M(u_0)$. Choose $R$ once and for all sufficiently large that
\[
CR^{-2s}\big(m+m^{N/(N-2s)}\big)\le4se,
\]
where $C$ is the constant in Proposition \ref{prop:coercive-virial}. Then
\begin{equation}
\mathcal M_R'(t)\le-4se-c_{N,s} \mathcal D_R(u(t)).
\label{eq:negative-MR-prime}
\end{equation}

Assume for contradiction that $T_+=\infty$. Integrating the constant negative term in \eqref{eq:negative-MR-prime},
\[
 \mathcal M_R(t)\le  \mathcal M_R(0)-4set.
\]
Hence $ \mathcal M_R(t)\to-\infty$ as $t\to\infty$. Set
\[
y(t):=- \mathcal M_R(t).
\]
There exists $t_0>0$ such that
\begin{equation}
y(t)\ge2C_0m,
\quad t\ge t_0,
\label{eq:y-large}
\end{equation}
where $C_0$ is the first constant in Proposition \ref{prop:virial-defect}. By \eqref{eq:virial-defect} and \eqref{eq:y-large},
\[
\frac12y(t) \le C_1R m^{1-1/(2s)} \mathcal D_R(u(t))^{1/(2s)}.
\]
Consequently,
\begin{equation}
\mathcal D_R(u(t)) \ge cR^{-2s}m^{-(2s-1)}y(t)^{2s}.
\label{eq:D-lower-y}
\end{equation}
On the other hand, \eqref{eq:negative-MR-prime} implies
\begin{equation}
y'(t)=- \mathcal M_R'(t)\ge c_0 \mathcal D_R(u(t)).
\label{eq:y-prime-D}
\end{equation}
Combining \eqref{eq:D-lower-y} and \eqref{eq:y-prime-D},
\[
y'(t) \ge cR^{-2s}m^{-(2s-1)}y(t)^{2s}, \quad t\ge t_0.
\]
Since $2s>1$, integration gives
\[
y(t)^{1-2s}\le y(t_0)^{1-2s}
-c(2s-1)R^{-2s}m^{-(2s-1)}(t-t_0).
\]
For a sufficiently large finite $t$, the right-hand side is nonpositive. The left-hand side, however, is positive whenever $u(t)\in H^s(\R^N)$, since the virial form is finite. Hence $T_+<\infty$.

For the backward direction, define $\widetilde u(t,x)=\overline{u(-t,x)}$. Since the equation is invariant under time reversal combined with complex conjugation, $\widetilde u$ solves \eqref{eq:fnls}, is radial, and has the same mass and energy as $u$. Applying the forward argument to $\widetilde u$ gives $T_-<\infty$.
\end{proof}

\begin{proof}[Proof of Corollary~\ref{cor:energy-space}]
Let $E(u_0)=-e<0$ and $m=M(u_0)>0$. Fix $R$ sufficiently large that the error coefficient in \eqref{eq:energy-virial-integrated} is at most $4se$. Proposition~\ref{prop:energy-virial}
then gives, on every compact subinterval of the lifespan,
\begin{equation}
\mathcal M_R(u(t_2))-\mathcal M_R(u(t_1))
\le -4se(t_2-t_1)-c\int_{t_1}^{t_2}\mathcal D_R(u(t)) dt.
\label{eq:energy-virial-negative}
\end{equation}
Suppose $T_+=\infty$ and set $y(t)=-\mathcal M_R(u(t))$. By \eqref{eq:energy-virial-negative}, $y(t)$ grows at least linearly, so there is
$t_0>0$ such that $y(t)\ge2C_0m$ for $t\ge t_0$, with $C_0$
as in Proposition~\ref{prop:virial-defect}. That proposition yields
\[
\mathcal D_R(u(t))\ge cR^{-2s}m^{-(2s-1)}y(t)^{2s}
\quad(t\ge t_0).
\]
Consequently, with $K=cR^{-2s}m^{-(2s-1)}>0$,
\[
y(t_2)-y(t_1)\ge K\int_{t_1}^{t_2}y(t)^{2s} dt,
\quad t_2\ge t_1\ge t_0.
\]
Set $z(t)=y(t_0)+K\int_{t_0}^t y(\tau)^{2s} d\tau$.
Then $y(t)\ge z(t)>0$ and $z'(t)\ge Kz(t)^{2s}$.
Since $2s>1$,
\[
z(t)^{1-2s}\le z(t_0)^{1-2s}-K(2s-1)(t-t_0).
\]
The right-hand side becomes nonpositive at a finite time, contrary to the positivity of $z$ on the lifespan. Thus $T_+<\infty$. Time reversal gives $T_-<\infty$, and the norm divergence follows from Proposition~\ref{prop:local-theory}.
\end{proof}

\begin{remark}[The half-wave endpoint]
At $s=1/2$, the final differential inequality would have the form $y'\gtrsim y$, which yields exponential growth but does not force a finite-time singularity. Thus the strict condition $s>1/2$ is reflected directly in the Osgood closure.
\end{remark}

\section{Strong instability of ground-state standing waves}
\label{sec:instability}

Let $Q_\omega$ be a positive radial ground state of \eqref{eq:ground-state}. At the mass-critical exponent, the ground state has zero energy. In fact, multiplying \eqref{eq:ground-state} by $Q_\omega$ and integrating over $\R^N$, we obtain the Nerahi identity
\[
\|(-\Delta)^{s/2}Q_\omega\|_{L^2}^2
+\omega\|Q_\omega\|_{L^2}^2
=\|Q_\omega\|_{L^q}^q.
\]
Moreover, the fractional Pohozaev identity for \eqref{eq:ground-state} gives
\[
\frac{N-2s}{2}
\|(-\Delta)^{s/2}Q_\omega\|_{L^2}^2
+\frac{N\omega}{2}\|Q_\omega\|_{L^2}^2
=\frac{N}{q}\|Q_\omega\|_{L^q}^q.
\]
This identity is justified by the standard regularity and decay properties of fractional ground states; see, for example, \cite{FrankLenzmannSilvestre}.

Subtracting $N/2$ times the Nehari identity from the Pohozaev identity, we obtain
\[
-s\|(-\Delta)^{s/2}Q_\omega\|_{L^2}^2
=-\frac{N(q-2)}{2q}\|Q_\omega\|_{L^q}^q.
\]
Since $N(q-2)=4s$, it follows that
\[
\frac12\|(-\Delta)^{s/2}Q_\omega\|_{L^2}^2
=\frac1q\|Q_\omega\|_{L^q}^q.
\]
Consequently,
\[
E(Q_\omega)=0
\]
hence
\begin{equation}
\frac12\|(-\Delta)^{s/2}Q_\omega\|_{L^2}^2 =\frac1q\|Q_\omega\|_{L^q}^q.
\label{eq:ground-pohozaev}
\end{equation}

\begin{proof}[Proof of Corollary \ref{cor:instability}]
Fix $c>1$ and set
\[
u_{0,c}:=cQ_\omega.
\]
Then
\[
\|u_{0,c}-Q_\omega\|_{H^s} =(c-1)\|Q_\omega\|_{H^s}\longrightarrow0
\quad\text{as }c\downarrow1.
\]
By \eqref{eq:ground-pohozaev},
\begin{align*}
E(cQ_\omega)
&=\frac{c^2}{2}\|(-\Delta)^{s/2}Q_\omega\|_{L^2}^2 -\frac{c^q}{q}\|Q_\omega\|_{L^q}^q\\
&=\frac12(c^2-c^q)\|(-\Delta)^{s/2}Q_\omega\|_{L^2}^2<0,
\end{align*}
since $q>2$ and $c>1$. The datum $cQ_\omega$ is radial, so Corollary \ref{cor:energy-space} implies finite-time blow-up. Choosing $c>1$ sufficiently close to $1$ gives the required perturbation inside any prescribed $H^s(\R^N)$ neighborhood of $Q_\omega$.
\end{proof}

\appendix

\section{Construction of the dyadic cutoffs}
\label{app:cutoffs}

For completeness, we record a standard construction of the cutoff families used in Section~\ref{sec:defect}. Choose a nonnegative function
\[
\chi\in C_c^\infty((1/2,4))
\]
such that
\[
\chi(r)>0 \quad\text{for }1\le r\le2.
\]
For $\rho_j=2^jR$, define
\[
\widetilde\zeta_j(x) := \chi\left(\frac{|x|}{\rho_j}\right).
\]
For every $x\ne0$ there exists $j$ such that
$1\le |x|/\rho_j\le2$, so
\[
S(x):= \left(\sum_{k\in\mathbb Z}\widetilde\zeta_k(x)^2\right)^{1/2}
\]
is strictly positive. Since $\chi$ has compact support, only a
uniformly bounded number of terms contribute at each $x\ne0$.
Consequently,
\[
0<c\le S(x)\le C
\]
uniformly in $x\ne0$ and $R>0$. Set
\[
\zeta_j(x) := \frac{\widetilde\zeta_j(x)}{S(x)}.
\]
Then
\[
\sum_{j\in\mathbb Z}\zeta_j(x)^2=1, \quad x\ne0,
\]
and
\[
\supp\zeta_j \subset \left\{x:\frac12\rho_j<|x|<4\rho_j\right\}.
\]
By dyadic scaling and the uniform bounds on $S$ and its derivatives,
\[
|\nabla^k\zeta_j(x)| \le C_k\rho_j^{-k}, \quad k=1,2.
\]
In particular, the supports $\{\supp\zeta_j\}$ have uniformly finite
overlap. 

To construct the enlarged cutoffs, choose
\[
\eta\in C_c^\infty((1/4,8)), \quad 0\le\eta\le1,
\]
with
\[
\eta(r)=1 \quad\text{for }r\in[1/2,4],
\]
and define
\[
\eta_j(x):=\eta\left(\frac{|x|}{\rho_j}\right).
\]
Then
\[
\eta_j\equiv1 \quad\text{on }\supp\zeta_j,
\]
while
\[
\supp\eta_j \subset \left\{x \in \R^N:\frac14\rho_j<|x|<8\rho_j\right\},
\]
and
\[
|\nabla^k\eta_j(x)| \le C_k\rho_j^{-k}, \quad k=1,2.
\]
The supports $\{\supp\eta_j\}$ also have uniformly finite overlap. Thus the two cutoff families satisfy all properties required in
Section~\ref{sec:defect}.

\section{Proof of Lemma \ref{lem:first-order-form}}
\label{app:first-order-form}

The argument refines the proof of \cite[Lemma~A.1]{BHL}. We first take $f\in\mathcal S(\R^N)$ and set
\[
A:=\|a\|_{L^\infty}, \quad L:=\|\nabla a\|_{L^\infty}.
\]
The case $A=0$ is trivial, so assume $A>0$.

\medskip
\noindent\emph{Step 1: Reduction to a fractional product estimate.}
For each component $a_\ell$,
\[
\int_{\R^N}\overline f a_\ell\partial_\ell f dx
=
\left\langle (-\Delta)^{1/4}(a_\ell f), (-\Delta)^{-1/4}\partial_\ell f \right\rangle.
\]
Since $(-\Delta)^{-1/4}\partial_\ell=R_\ell (-\Delta)^{1/4}$ and the Riesz transform $R_\ell$ is bounded on $L^2(\R^N)$,
\[
\left|
\int_{\R^N}\overline f a\cdot\nabla f dx
\right|
\lesssim_N \sum_{\ell=1}^N \|(-\Delta)^{1/4}(a_\ell f)\|_{L^2} \|(-\Delta)^{1/4}f\|_{L^2}.
\]
It therefore suffices to prove, for a scalar bounded Lipschitz coefficient $b$ satisfying $\|b\|_{L^\infty}\le A$ and
$\|\nabla b\|_{L^\infty}\le L$,
\begin{equation}
\|(-\Delta)^{1/4}(bf)\|_{L^2} \lesssim A\|(-\Delta)^{1/4}f\|_{L^2}
+ (AL)^{1/2}\|f\|_{L^2}.
\label{eq:half-product}
\end{equation}

\medskip
\noindent\emph{Step 2: The $H^{1/2}$ product estimate.}
Using the Gagliardo characterization of the homogeneous $H^{1/2}$ seminorm,
\[
\|(-\Delta)^{1/4}g\|_{L^2}^2 \simeq \iint_{\R^N\times\R^N} \frac{|g(x)-g(y)|^2}{|x-y|^{N+1}} dx dy,
\]
we write
\[
b(x)f(x)-b(y)f(y)=b(x)(f(x)-f(y))+(b(x)-b(y))f(y).
\]
Hence
\[
\begin{aligned}
\|(-\Delta)^{1/4}(bf)\|_{L^2}^2 \lesssim  A^2\|(-\Delta)^{1/4}f\|_{L^2}^2+\int_{\R^N}|f(y)|^2 \left(\int_{\R^N}\frac{|b(x)-b(y)|^2}{|x-y|^{N+1}} dx\right)dy.
\end{aligned}
\]
Boundedness and Lipschitz continuity give
\[
|b(x)-b(y)| \le \min\{2A,L|x-y|\}.
\]
Thus, for fixed $y$, polar coordinates yield
\[
\int_{\R^N} \frac{|b(x)-b(y)|^2}{|x-y|^{N+1}} dx
\lesssim_N \int_0^\infty \frac{\min\{4A^2,L^2r^2\}}{r^2} dr.
\]
If $L=0$ this integral vanishes. Otherwise, splitting at $r_0=2A/L$ gives
\[
L^2\int_0^{r_0}dr+4A^2\int_{r_0}^\infty r^{-2} dr=4AL.
\]
Therefore
\[
\|(-\Delta)^{1/4}(bf)\|_{L^2}^2 \lesssim A^2\|(-\Delta)^{1/4}f\|_{L^2}^2
+ AL\|f\|_{L^2}^2,
\]
which implies \eqref{eq:half-product}.

\medskip
\noindent\emph{Step 3: Conclusion.}
Combining Steps 1 and 2,
\[
\left|\int\overline f a\cdot\nabla f \right|
\lesssim
A\|(-\Delta)^{1/4}f\|_{L^2}^2+(AL)^{1/2}\|f\|_{L^2}\|(-\Delta)^{1/4}f\|_{L^2}.
\]
Young's inequality, applied as
\[
(AL)^{1/2}\|f\|_{L^2}\|(-\Delta)^{1/4}f\|_{L^2}
\le \frac12A\|(-\Delta)^{1/4}f\|_{L^2}^2+\frac12L\|f\|_{L^2}^2,
\]
gives
\[
\left|\int\overline f a\cdot\nabla f\right|
\lesssim A\|(-\Delta)^{1/4}f\|_{L^2}^2 +L\|f\|_{L^2}^2.
\]
Taking the imaginary part and absorbing the factor $2$ into the constant proves \eqref{eq:first-order-form} for Schwartz functions. The general $H^{1/2}$ case follows by density argument as explained in Remark \ref{rem:density-H1/2}. \hfill $\square$

\section{The virial identity at energy regularity}
\label{app:energy-approx}

The virial identity established for strong solutions involves derivatives that are not directly available for solutions in the energy space $H^s(\R^N)$. Although one might try to justify the identity by approximating the initial data with more regular data, this approach would require a suitable persistence-of-regularity result. Such persistence is not part of the $H^s$ local theory used here and may require additional care because the nonlinearity $F(z)=|z|^{4s/N}z$ has limited smoothness at the origin. We therefore avoid any assumption that an energy-space solution belongs to $H^{2s}(\R^N)$.

Instead, we regularize the solution in space, derive a forced virial identity for the regularized equation, and show that the forcing error vanishes using the lossless radial Strichartz bounds in Proposition~\ref{prop:local-theory}. We first establish the integrated identity for sufficiently smooth weights and then apply it to regularizations of the exponential multiplier.

Recall from \eqref{eq:p-q} and \eqref{eq:gamma-rho} that
\[
p=\frac{4s}{N},\quad q=p+2,\quad
\gamma=\frac{2(N+2s)}{N-2s},\quad
\rho=\frac{2N(N+2s)}{N^2+4s^2}.
\]
Set
\[
\nu=\frac{N\rho}{N-s\rho}=\frac{2(N+2s)}{N-2s}.
\]
When $N\ge2$ and $N/(2N-1)\le s<1$, we have $1<\rho'<2$, $0<s<N/2$, and
\begin{equation}
\frac1{\rho'}=\frac1\rho+\frac p\nu,
\quad
\frac{p+2}{\gamma}=1-\frac{2s}{N}<1.
\label{eq:energy-comm-exponents}
\end{equation}
For every compact interval $I$ contained in the lifespan, the Sobolev embedding $W^{s,\rho}\hookrightarrow L^\nu$ and Proposition~\ref{prop:local-theory} imply
\begin{equation}
\|u\|_{L^\gamma(I;L^\nu)}
\lesssim\|u\|_{L^\gamma(I;W^{s,\rho})}<\infty.
\label{eq:energy-comm-sobolev}
\end{equation}

\begin{lemma}[Spatial commutator bound]
\label{lem:spatial-virial-commutator}
Let $u\in L^\gamma(I;W^{s,\rho}(\R^N))$ and let $F(z)=|z|^pz$. Choose a real-valued radial function $\chi\in C_c^\infty(\R^N)$ such that $\chi=1$ near the origin. For $0<h\le1$, let $J_h$ denote the Fourier multiplier with symbol $\chi(h\xi)$, and define
\[
v_h=J_hu,\quad e_h=J_hF(u)-F(v_h).
\]
Then, for almost every $t\in I$,
\begin{align}
\|e_h(t)\|_{L^{\rho'}}
&\le Ch^s\|u(t)\|_{L^\nu}^p
       \|u(t)\|_{W^{s,\rho}},
\label{eq:energy-comm-pointwise}\\
\|\nabla v_h(t)\|_{L^\rho}
&\le Ch^{s-1}\|u(t)\|_{W^{s,\rho}},
\quad
\|v_h(t)\|_{L^\rho}\le C\|u(t)\|_{L^\rho}.
\label{eq:energy-comm-grad}
\end{align}
Consequently,
\begin{equation}
\begin{aligned}
\int_I\|e_h(t)\|_{L^{\rho'}}
  \bigl(\|\nabla v_h(t)\|_{L^\rho}
       +\|v_h(t)\|_{L^\rho}\bigr) dt
\le
C\bigl(h^{2s-1}+h^s\bigr)|I|^{2s/N}
\|u\|_{L^\gamma(I;W^{s,\rho})}^{p+2}.
\end{aligned}
\label{eq:energy-comm-integral}
\end{equation}
\end{lemma}

\begin{proof}
The convolution kernel of $J_h$ has an $L^1$ norm bounded independently of $h$, so $J_h$ is uniformly bounded on every $L^a$, $1\le a\le\infty$. Since $\chi=1$ near the origin, the Mihlin multiplier theorem also gives, for $1<a<\infty$ and $0<s<1$,
\[
\|(1-J_h)w\|_{L^a}
\lesssim h^s\|(-\Delta)^{s/2} w\|_{L^a},\quad
\|\nabla J_hw\|_{L^a}
\lesssim h^{s-1}\|(-\Delta)^{s/2}w\|_{L^a}.
\]

The inequality
\[
|F(z)-F(w)|\lesssim (|z|^p+|w|^p)|z-w|
\]
and the first identity in \eqref{eq:energy-comm-exponents} yield
\[
\begin{aligned}
\|F(u)-F(v_h)\|_{L^{\rho'}}
&\lesssim \bigl(\|u\|_{L^\nu}^p + \|v_h\|_{L^\nu}^p\bigr)\|u-v_h\|_{L^\rho}\\
&\lesssim h^s\|u\|_{L^\nu}^p \|(-\Delta)^{s/2}u\|_{L^\rho}.
\end{aligned}
\]
Moreover, the fractional chain rule for $0<s<1$ (see e.g. \cite[Lemma~2.1]{DinhMassCritical}) gives
\[
\begin{aligned}
\|(-\Delta)^{s/2}F(v_h)\|_{L^{\rho'}}
\lesssim\|v_h\|_{L^\nu}^p \|(-\Delta)^{s/2}v_h\|_{L^\rho}
\lesssim\|u\|_{L^\nu}^p \|(-\Delta)^{s/2} u\|_{L^\rho}.
\end{aligned}
\]
Since
\[
e_h=J_h\bigl(F(u)-F(v_h)\bigr)+(J_h-1)F(v_h),
\]
these estimates prove \eqref{eq:energy-comm-pointwise}. The multiplier bounds also imply \eqref{eq:energy-comm-grad}.

Finally, H\"older's inequality in time and \eqref{eq:energy-comm-exponents} give
\[
\begin{aligned}
\int_I\|u(t)\|_{L^\nu}^p  \|u(t)\|_{W^{s,\rho}}^2 dt
\le |I|^{2s/N} \|u\|_{L^\gamma(I;L^\nu)}^p \|u\|_{L^\gamma(I;W^{s,\rho})}^2.
\end{aligned}
\]
Combining this with \eqref{eq:energy-comm-sobolev} proves \eqref{eq:energy-comm-integral}.
\end{proof}

We next establish the integrated identity for weights whose derivatives are sufficiently regular to justify the classical commutator calculation. 

\begin{lemma}[Integrated virial identity for smooth weights]
\label{lem:Hs-smooth-virial}
Assume \eqref{eq:expo-energy} and let $u$ be the radial energy-space solution given by Proposition~\ref{prop:local-theory}. Let $\phi\in C^4(\R^N;\R)$ satisfy $\nabla\phi\in W^{3,\infty}(\R^N)$ and
\[
|\Delta^2\phi(x)|\le C_\phi|x|^{-2(1-s)},
\quad x\ne0.
\]
Then, for every $-T_-<t_1\le t_2<T_+$,
\begin{align}
\mathcal M_\phi(u(t_2))-\mathcal M_\phi(u(t_1))
={}&4\int_{t_1}^{t_2} \int_0^\infty \mu^s
\int_{\R^N}\partial_k\overline{u_\mu} 
\partial_{kl}\phi \partial_lu_\mu dx d\mu dt
\notag\\
&-\int_{t_1}^{t_2} \int_0^\infty \mu^s
\int_{\R^N}(\Delta^2\phi)|u_\mu|^2 dx d\mu dt
\notag\\
&-\frac{2p}{q}\int_{t_1}^{t_2}
\int_{\R^N}(\Delta\phi)|u|^q dx dt.
\label{eq:energy-virial-identity}
\end{align}
Here and below, repeated spatial indices are summed. The virial form on the left is understood through its continuous extension to $H^s(\R^N)$, as in Lemma~\ref{lem:first-order-form}.
\end{lemma}

\begin{proof}
Fix a compact interval $I=[t_1,t_2]$ within the lifespan and put $v_h=J_hu$, as in Lemma~\ref{lem:spatial-virial-commutator}. Since $J_h$ is uniformly bounded on $H^s(\R^N)$, converges strongly to the identity, and $u(I)$ is compact in $H^s$, we have
\[
v_h\longrightarrow u\quad\text{in }C(I;H^s(\R^N)) \quad \text{as } h \downarrow 0.
\]

We first justify differentiation of the virial action for $v_h$. The embedding $H^s(\R^N)\hookrightarrow L^q(\R^N)$ and H\"older's inequality imply
\[
\begin{aligned}
\|F(f)-F(g)\|_{H^{-s}}
&\lesssim\|F(f)-F(g)\|_{L^{q'}}\\
&\lesssim\bigl(\|f\|_{L^q}^p +\|g\|_{L^q}^p\bigr) \|f-g\|_{L^q}.
\end{aligned}
\]
It follows that $F(u)\in C(I;H^{-s}(\R^N))$. Since $(-\Delta)^su\in C(I;H^{-s}(\R^N))$, the equation in its Duhamel form gives $u\in C^1(I;H^{-s}(\R^N))$. For each fixed $h>0$, the compactly supported Fourier multiplier $J_h$ maps $H^{-s}(\R^N)$ into $H^k(\R^N)$ for every $k\ge0$. Consequently,
\[
v_h\in C^1(I;H^k(\R^N))\quad\text{for every }k\ge0.
\]
The regularized equation is
\[
i\partial_tv_h=(-\Delta)^sv_h-F(v_h)-e_h,
\quad e_h=J_hF(u)-F(v_h).
\]
For fixed $h$, the function $v_h$ belongs to $H^k(\R^N)$ for every $k$ and hence to $H^1(\R^N)\cap L^\infty(\R^N)$. Although $F$ need not be twice differentiable at the origin, it is $C^1$ and $|DF(z)|\lesssim|z|^p$. Thus
\[
\|F(v_h)\|_{L^2}+\|\nabla F(v_h)\|_{L^2}
\lesssim\|v_h\|_{L^\infty}^p
\bigl(\|v_h\|_{L^2}+\|\nabla v_h\|_{L^2}\bigr).
\]
Since $J_hF(u)$ belongs to every $H^k(\R^N)$, we have $e_h\in H^1(\R^N)$. This regularity, together with the boundedness of $\nabla\phi$ and $\Delta\phi$, justifies the differentiations and integrations by parts below. If necessary, the latter may first be performed with a spatial cutoff and then the cutoff removed using the corresponding $L^2$ bounds.

Write $L=-\Delta$ and introduce the formally skew-adjoint first-order operator
\[
\mathcal C_\phi :=\nabla\phi\cdot\nabla+\frac12\Delta\phi.
\]
With the $L^2$ inner product taken linear in its second argument, the virial action satisfies
\[
\mathcal M_\phi(v)=2\langle v,-i\mathcal C_\phi v\rangle
\]
for sufficiently regular $v$. The linear part of the equation contributes
$2\langle v_h,[L^s,\mathcal C_\phi]v_h\rangle$ to its derivative.

To justify the fractional commutator, first truncate the Balakrishnan representation by setting
\[
L^s_{\varepsilon,M}:=c_s^2\int_\varepsilon^M \mu^{s-1}L(L+\mu)^{-1} d\mu,
\quad 0<\varepsilon<M<\infty.
\]
The resolvent identity yields, on sufficiently regular functions,
\[
[L^s_{\varepsilon,M},\mathcal C_\phi]
=c_s^2\int_\varepsilon^M\mu^s
(L+\mu)^{-1}[L,\mathcal C_\phi]
(L+\mu)^{-1} d\mu.
\]
Direct differentiation gives
\[
[L,\mathcal C_\phi]f
=-2\partial_{jk}\phi \partial_{jk}f
-2\nabla\Delta\phi\cdot\nabla f
-\frac12(\Delta^2\phi)f.
\]
Integrating by parts, we obtain
\[
\langle f,[L,\mathcal C_\phi]f\rangle
=2\int_{\R^N}\partial_j\overline f 
  \partial_{jk}\phi \partial_kf dx
-\frac12\int_{\R^N}(\Delta^2\phi)|f|^2 dx.
\]
With $v_{h,\mu}=c_s(L+\mu)^{-1}v_h$, the truncated commutator therefore contributes
\[
\begin{aligned}
4\int_\varepsilon^M\mu^s 
  \int_{\R^N}\partial_j\overline{v_{h,\mu}} 
  \partial_{jk}\phi \partial_kv_{h,\mu} dx d\mu -\int_\varepsilon^M\mu^s 
  \int_{\R^N}(\Delta^2\phi)|v_{h,\mu}|^2 dx d\mu.
\end{aligned}
\]
The first integral converges absolutely as $\varepsilon\downarrow0$ and $M\uparrow\infty$, because $D^2\phi\in L^\infty(\R^N)$ and the Balakrishnan identity gives
\[
\int_0^\infty\mu^s\|\nabla v_{h,\mu}\|_{L^2}^2 d\mu
=s\|v_h\|_{\dot H^s}^2.
\]
For the lower-order integral, the assumed bound on $\Delta^2\phi$ and the fractional Hardy inequality with exponent $1-s$ give (see \eqref{eq:iden-L2})
\[
\begin{aligned}
\int_0^\infty\mu^s 
  \int_{\R^N}|\Delta^2\phi| |v_{h,\mu}|^2 dx d\mu
  \lesssim C_\phi
  \int_0^\infty\mu^s
  \|(-\Delta)^{(1-s)/2}v_{h,\mu}\|_2^2 d\mu
=C_\phi s\|v_h\|_2^2.
\end{aligned}
\]
It remains to justify convergence of the truncated commutator to $[L^s,\mathcal C_\phi]$. For fixed $h>0$, we have $v_h\in H^k(\R^N)$ for every $k\ge0$. Moreover, since $\nabla\phi\in W^{3,\infty}(\R^N)$, the definition of $\mathcal C_\phi$ implies that
\[
\mathcal C_\phi v_h\in H^2(\R^N)
\subset H^{2s}(\R^N).
\]
Thus both $v_h$ and $\mathcal C_\phi v_h$ belong to the domain of $L^s$.

To justify the convergence, define
\[
m_{\varepsilon,M}(\lambda)
:=c_s^2\int_\varepsilon^M
\mu^{s-1}\frac{\lambda}{\lambda+\mu}\,d\mu,
\qquad \lambda\ge0.
\]
By the Balakrishnan formula,
\[
0\le m_{\varepsilon,M}(\lambda)
\le\lambda^s,
\qquad
m_{\varepsilon,M}(\lambda)\longrightarrow\lambda^s
\]
as $\varepsilon\downarrow0$ and $M\uparrow\infty$. The spectral theorem and dominated convergence therefore give
\[
L^s_{\varepsilon,M}w\longrightarrow L^sw
\quad\text{in }L^2(\R^N)
\]
for every $w\in H^{2s}(\R^N)$.

Applying this convergence to $v_h$ and $\mathcal C_\phi v_h$, and using the formal skew-adjointness of $\mathcal C_\phi$, we obtain
\[
\big\langle v_h,
[L^s_{\varepsilon,M},\mathcal C_\phi]v_h
\big\rangle
\longrightarrow
\big\langle v_h,
[L^s,\mathcal C_\phi]v_h
\big\rangle.
\]
Since the two resolvent quadratic forms have already been shown to converge absolutely, we may remove the cutoffs and recover the
fractional kinetic terms in \eqref{eq:energy-virial-identity}.

For the nonlinear part, the $C^1$ regularity of $F$ gives the pointwise identity
\[
\begin{aligned}
2\operatorname{Re}\bigl[
  -\overline{F(v_h)} \nabla\phi\cdot\nabla v_h
  +\overline{v_h} \nabla\phi\cdot\nabla F(v_h)
  \bigr] =\frac{2p}{q}\nabla\phi\cdot\nabla(|v_h|^q).
\end{aligned}
\]
This identity holds also where $v_h=0$, by continuity, and does not require differentiating the potentially nonsmooth function $r\mapsto r^{p/2}$ at zero. Since $v_h\in H^1(\R^N)\cap L^\infty(\R^N)$, we have $|v_h|^q\in W^{1,1}(\R^N)$. Integration by parts consequently yields
\[
-\frac{2p}{q}\int_{\R^N}(\Delta\phi)|v_h|^q dx.
\]

The forcing term in the regularized equation contributes
\[
\begin{aligned}
\mathcal E_h(t)
&=2\operatorname{Im}\int_{\R^N} 
\bigl[-i\overline{e_h} \nabla\phi\cdot\nabla v_h
+i\overline{v_h} \nabla\phi\cdot\nabla e_h\bigr]dx\\
&=-4\operatorname{Re}\int_{\R^N}
\overline{e_h} \nabla\phi\cdot\nabla v_h dx
-2\operatorname{Re}\int_{\R^N}
(\Delta\phi)\overline{v_h}e_h dx.
\end{aligned}
\]
Combining the preceding computations gives the forced integrated virial identity for $v_h$, with the right-hand side of \eqref{eq:energy-virial-identity} evaluated at $v_h$ and the additional term $\int_I\mathcal E_h(t) dt$. By H\"older's inequality and Lemma~\ref{lem:spatial-virial-commutator},
\[
\begin{aligned}
\int_I|\mathcal E_h(t)| dt
&\lesssim
(\|\nabla\phi\|_\infty+\|\Delta\phi\|_\infty)
\int_I\|e_h(t)\|_{L^{\rho'}}
\bigl(\|\nabla v_h(t)\|_{L^\rho}
     +\|v_h(t)\|_{L^\rho}\bigr) dt\\
&\lesssim_\phi
\bigl(h^{2s-1}+h^s\bigr)|I|^{2s/N}
\|u\|_{L^\gamma(I;W^{s,\rho})}^{p+2}
\longrightarrow0,
\end{aligned}
\]
because $s>1/2$.

It remains to pass to the limit in the unforced terms. Denote the principal Balakrishnan quadratic form by $\mathcal K_\phi$. By polarization and \eqref{eq:balakrishnan-energy},
\[
\begin{aligned}
|\mathcal K_\phi(v_h)-\mathcal K_\phi(u)|
\lesssim{}&\|D^2\phi\|_{L^\infty}
\bigl(\|(-\Delta)^{s/2}v_h\|_{L^2}+\|(-\Delta)^{s/2}u\|_{L^2}\bigr)
\|(-\Delta)^{s/2}(v_h-u)\|_{L^2}.
\end{aligned}
\]
For the lower-order quadratic form $\mathcal L_\phi$, the preceding Hardy--Balakrishnan estimate and polarization imply
\[
|\mathcal L_\phi(v_h)-\mathcal L_\phi(u)|
\lesssim C_\phi
\bigl(\|v_h\|_{L^2}+\|u\|_{L^2}\bigr)\|v_h-u\|_{L^2}.
\]
Both differences tend to zero uniformly for $t\in I$ because $v_h\to u$ in $C(I;H^s(\R^N))$. The nonlinear term converges as well, since $v_h\to u$ in $C(I;L^q(\R^N))$, which implies $|v_h|^q\to|u|^q$ in $C(I;L^1(\R^N))$, and $\Delta\phi\in L^\infty(\R^N)$. Finally, Lemma~\ref{lem:first-order-form} yields convergence of the endpoint virial actions. We may therefore let $h\downarrow0$ in the integrated identity and obtain \eqref{eq:energy-virial-identity}.
\end{proof}

\begin{proof}[Proof of Proposition~\ref{prop:energy-virial}]
Fix $R>0$ and a compact interval $[t_1,t_2]$ in the lifespan. We first apply Lemma~\ref{lem:Hs-smooth-virial} with the regularized multiplier $\phi_{R,\delta}$ from Section~\ref{subsec:regularization}, where $0<\delta<R$. For each fixed $\delta$, this multiplier is smooth and $\nabla\phi_{R,\delta}\in W^{3,\infty}(\R^N)$. Moreover, \eqref{eq:regularized-bounds2} yields, uniformly in $\delta$,
\[
|\Delta^2\phi_{R,\delta}(x)|
\lesssim R^{-2s}(|x|+\delta)^{2s-2}
\le R^{-2s}|x|^{-2(1-s)},\quad x\ne0.
\]
Thus all the hypotheses of the lemma are satisfied. Since the solution-regularization limit $h\downarrow0$ has already been taken there, we now let $\delta\downarrow0$.

The Hessians $D^2\phi_{R,\delta}$ converge pointwise to $D^2\phi_R$ and remain uniformly bounded. The principal Balakrishnan integrands are therefore dominated by a constant multiple of $\mu^s|\nabla u_\mu|^2$, which is integrable by \eqref{eq:balakrishnan-energy} and $u\in C([t_1,t_2];H^s(\R^N))$. Similarly, $\Delta\phi_{R,\delta}\to\Delta\phi_R$ pointwise, and the uniform Laplacian bound permits dominated convergence in the nonlinear term.

For the lower-order Balakrishnan term, the uniform bilaplacian estimate and the fractional Hardy inequality give
\[
\begin{aligned}
\int_{t_1}^{t_2} \int_0^\infty\mu^s
\int_{\R^N}|x|^{-2(1-s)}|u_\mu(t,x)|^2
 dx d\mu dt 
\lesssim (t_2-t_1)
\sup_{t\in[t_1,t_2]}\|u(t)\|_{L^2}^2<\infty.
\end{aligned}
\]
This furnishes an integrable majorant independent of $\delta$, so dominated convergence also applies to the lower-order term.

It remains to justify the endpoint virial actions. The gradients $\nabla\phi_{R,\delta}$ converge pointwise to $\nabla\phi_R$, while their $L^\infty$ norms and the $L^\infty$ norms of their first derivatives are bounded uniformly in $\delta$. For a Schwartz function, convergence of the corresponding virial actions follows from dominated convergence. Lemma~\ref{lem:first-order-form} gives a uniform quadratic-form bound on $H^{1/2}(\R^N)$, and a density argument extends the convergence to $u(t_1),u(t_2)\in H^s(\R^N)$. We conclude that \eqref{eq:energy-virial-identity} remains valid in integrated form with $\phi=\phi_R$.

We now apply the algebraic estimates established in the proof of Proposition~\ref{prop:coercive-virial}. For every fixed time, Lemma~\ref{lem:geometry}, Proposition~\ref{prop:defect-transfer}, and Proposition~\ref{prop:weighted-interpolation} apply to the radial function $u(t)\in H^s(\R^N)$. Together with the lower-order Balakrishnan estimate, they bound the integrand of the energy-space virial identity by
\[
8sE(u(t))-c_{N,s}\mathcal D_R(u(t)) +C_{N,s}R^{-2s} \left[M(u(t))+M(u(t))^{N/(N-2s)}\right].
\]
This step uses only the algebraic inequalities for the virial integrand; it does not invoke the strong-solution statement of Proposition~\ref{prop:coercive-virial}.

To justify integration of the defect term, Proposition~\ref{prop:defect-transfer} and \eqref{eq:balakrishnan-energy} give, for each fixed $R>0$,
\[
\mathcal D_R(v)
\lesssim\mathcal Q_R(v)+R^{-2s}\|v\|_{L^2}^2
\lesssim\|(-\Delta)^{s/2}v\|_{L^2}^2+R^{-2s}\|v\|_{L^2}^2
\lesssim_R\|v\|_{H^s}^2.
\]
Since $u\in C([t_1,t_2];H^s(\R^N))$, it follows that
\[
\int_{t_1}^{t_2}\mathcal D_R(u(t)) dt
\lesssim_R(t_2-t_1)
\sup_{t\in[t_1,t_2]}\|u(t)\|_{H^s}^2<\infty.
\]
We may therefore integrate the preceding virial estimate over $[t_1,t_2]$. Conservation of mass and energy then yields \eqref{eq:energy-virial-integrated}, as required.
\end{proof}

\section*{Declaration of AI use}

The authors used ChatGPT (GPT-5.6 Sol, OpenAI) during the preparation of this manuscript for language editing and discussions of possible approaches to certain technical arguments. The authors take full responsibility for the mathematical correctness and content of the paper.

\section*{Data availability statement}

Data sharing is not applicable to this article as no new data were created or analyzed in this study.

\end{document}